\documentclass[twoside,12pt]{article}
\usepackage{amssymb,amsmath,bm,euscript}
\usepackage{graphicx,color,caption}
\usepackage{ifpdf}

\usepackage[colorlinks=true]{hyperref}

\usepackage{braket}

\usepackage{algorithm, algorithmic}

\newtheorem{proposition}{Proposition}[section]
\newtheorem{theorem}[proposition]{Theorem}
\newtheorem{lemma}[proposition]{Lemma}

\newtheorem{definition}[proposition]{Definition}

\renewcommand{\theequation}{\thesection.\arabic{equation}}

\def\0{{\bf 0}}

\title{\bf{An Efficient Two-Sided Sketching Method for Large-Scale Tensor Decomposition Based on  Transformed Domains\thanks{\tiny Submitted to the editors DATE.
			
			\textbf{Funding}:{This work was supported in part by National Natural Science Foundation of China (No. 12071104) and Natural Science Foundation of Zhejiang Province. (No. LD19A010002).}}}}

\author{ \hspace{1mm}Zhiguang Cheng\thanks{\tiny Department of Mathematics, Hangzhou Dianzi University, Hangzhou, 310018, China; Department of Mathematics, Quzhou University, Quzhou, 324000, China; ({\tt
			32086@qzc.edu.cn}).}
	\and
	Gaohang Yu\thanks{\tiny Department of Mathematics, Hangzhou Dianzi University, Hangzhou, 310018, China; ({\tt
			maghyu@163.com}). }
		\and
	Xiaohao Cai \thanks{\tiny School of Electronics and Computer Science, University of Southampton, Southampton, UK; ({\tt
			x.cai@soton.ac.uk}). }
		\and
		Liqun Qi \thanks{\tiny Department of Applied Mathematics, The Hong Kong Polytechnic University, Kowloon, Hong Kong; ({\tt
				liqun.qi@polyu.edu.hk}). }
}

\begin{document}
\maketitle

\begin{abstract}
	Large tensors are frequently encountered in various fields such as computer vision, scientific simulations, sensor networks, and data mining. However, these tensors are often too large for convenient processing, transfer, or storage. Fortunately, they typically exhibit a low-rank structure that can be leveraged through tensor decomposition. However, performing large-scale tensor decomposition can be time-consuming. Sketching is a useful technique to reduce the dimensionality of the data. In this paper, we propose a novel two-sided sketching method based on the $\star_{L}$-product decomposition and transformed domains like the discrete cosine transformation. A rigorous theoretical analysis is also conducted to assess the approximation error of the proposed method. Specifically, we improve our method with power iteration to achieve more precise approximate solutions. Extensive numerical experiments and comparisons on low-rank approximation of synthetic large tensors and real-world data like color images and grayscale videos illustrate the efficiency and effectiveness of the proposed approach in terms of both CPU time and approximation accuracy.
\end{abstract}



\section{Introduction}
Multidimensional arrays, known as tensors, are often used to represent real-world high-dimensional data, such as videos \cite{T.-X17,T.-K07}, hyperspectral images \cite{Du16,Renard08,Wu22,Wang22}, multilinear signals \cite{Lathauwer97,Comon02}, and communication networks \cite{Papalexakis14,Nakatsuji17}. In most cases, these tensor data usually have a low-rank structure and can be approximated by tensor decomposition. Nevertheless, computing the tensor decomposition of these large-scale data is usually computationally demanding, and thus finding an accurate approximation of large-scale data with great efficiency plays a key role in tensor data analysis. Sketching is a useful technique for data compression, utilizing random projections or sampling to approximate the original data. Although the sketching technique may slightly reduce the accuracy of the approximation, it can significantly reduce the computational and storage complexity \cite{Liu22}. As a result, sketching is commonly used in low-rank tensor approximation, and many researchers have proposed various tensor sketching algorithms \cite{Cao23,Li17,WangY15,Che19,Che21,Che21J,Malik18,Ravishankar21,Sun20,Hur22,Qi21,Minster20,Dong23}. They have also been successfully applied to a variety of tasks, such as Kronecker product regression, polynomial approximation, the construction of deep convolutional neural networks \cite{Diao19,Han20,Kasiviswanathan18,Ma20,Shi19,Yang18}, etc.

Recently, extensive research has been carried out on the application of sketching algorithms for the low-rank matrix approximation \cite{Carmon19,David14,Tropp17,Tropp19}. Woodruff et al. \cite{David14} examined the numerical linear algebra algorithms of linear sketching techniques and identified their limitations. Tropp et al. \cite{Tropp17} developed a two-sided matrix sketching algorithm, which can maintain the structural properties of the input matrix and generate a low-rank matrix approximation with a given rank. Furthermore, Tropp et al. \cite{Tropp19} proposed a new matrix sketching algorithm to construct a low-rank approximation matrix from streaming data. These matrix sketching algorithms are very effective in reducing storage and computational costs when computing low-rank approximations of large-scale matrices. Subsequently, many researchers have applied matrix sketching algorithms to tensor decomposition and developed various low-rank tensor approximation algorithms. The followings are some low-rank tensor approximation algorithms based on different decompositions \cite{Cao23,Li17,WangY15,Che19,Che21,Che21J,Malik18,Ravishankar21,Sun20,Hur22,Qi21,Minster20,Dong23}.

Li et al. \cite{Li17} introduced a random algorithm for CANDECOMP/PARAFAC (CP) tensor decomposition in least-squares regression, aiming to achieve reduced dimensionality and sparsity in randomized linear mapping. Wang et al. \cite{WangY15} developed innovative techniques for performing randomized tensor contractions using the fast Fourier transform (FFT), avoiding explicit formation of tensors. The robust tensor power method based on the tensor sketch (TS-RTPM) can quickly explore the potential features of the tensor, but in some cases its approximation performance is limited. Cao et al. \cite{Cao23} proposed a data-driven framework called TS-RTPM-Net, which improves the accuracy of the estimation of TS-RTPM by jointly training the TS value matrices with the initial RTPM. 

Numerous studies have also been conducted based on Tucker decomposition \cite{Che19, Che21, Che21J, Malik18, Ravishankar21, Sun20, Minster20,Dong23}. For example, Che et al. \cite{Che19} devised an adaptive randomized approach to approximate Tucker decomposition. Malik et al. \cite{Malik18} proposed two randomized algorithms for low-rank Tucker decomposition, which entail a single pass of the input tensor by integrating sketching. Ravishankar et al. \cite{Ravishankar21} introduced the hybrid Tucker TensorSketch vector quantization (HTTSVQ) algorithm for dynamic light fields. Sun et al. \cite{Sun20} developed a randomized method for Tucker decomposition, which can provide a satisfactory approximation without a second pass on the original tensor data. Minster et al. \cite{Minster20} devised randomized adaptations of the THOSVD and STHOSVD algorithms. Dong et al. \cite{Dong23} presented two practical randomized algorithms for low-rank Tucker approximation of large tensors based on sketching and power scheme, with a rigorous error-bound analysis.

The work based on tensor-train (TT) decomposition can be found in \cite{Che19, Hur22, Kressner22, Che2022}. In particular, Che et al. \cite{Che19} designed an adaptive random algorithm to calculate the tensor column approximation. Hur et al. \cite{Hur22} introduced a sketching algorithm to construct a TT representation of a probability density from its samples, which can avoid the curse of dimensionality and sample complexities of the recovery problem. 
Qi and Yu \cite{Qi21} proposed a tensor sketching method based on the t-product, which can quickly obtain a low tubal rank tensor approximation. As pointed out by Kernfeld et al. \cite{Kernfeld15}, the t-product has a disadvantage in that, for real tensors, implementation of the t-product and factorizations using the t-product require intermediate complex arithmetic, which, even taking advantage of complex symmetry in the Fourier domain, is more expensive than real arithmetic. 


In this paper, based on the $\star_{L}$-product decomposition \cite{Kernfeld15},
we investigate two-sided sketching algorithms for low tubal rank tensor approximation. The main contributions of this paper are as follows. 
Firstly, we propose a new two-sided sketching method based on transformed domains, which can significantly improve the computational efficiency of the T-Sketch and rt-SVD methods, for low tubal rank approximation. Secondly, we establish a low tubal rank tensor approximation model based on the $\star_{L}$-product factorization, extending the two-sided matrix sketching algorithm proposed by Tropp et al. \cite{Tropp19} with subspace power iteration. Thirdly, a rigorous theoretical analysis is conducted to evaluate the approximation error of the proposed two-sided sketching method. 
Finally, extensive numerical experiments and comparisons on low-rank approximation of large tensors in different modelities (e.g. synthetic large tensors, color images,  and grayscale videos) illustrate the efficiency and effectiveness of the proposed approach in terms of both CPU time and approximation accuracy. 



The rest of this paper is organized as follows. Section \ref{sec:notation} introduces some common notations and preliminary. Our two-sided sketching algorithms (i.e., Algorithms \ref{alg:l-trp-sketch} and \ref{alg:DCT-Sketch-PI}) based on transformed domains are proposed in Section \ref{sec:our-ts-alg}. Section \ref{sec:err-ana} provides strict theoretical guarantees for the approximation error of the proposed algorithms. Section \ref{sec:experiments} presents detailed numerical experiments and comparisons, demonstrating the efficiency and effectiveness of the proposed algorithms. We conclude in Section \ref{sec:con}.  
Appendix A presents three matrix sketching techniques served as preliminary of the introduction of tensor sketching operators in Section \ref{subsec:tso}. Further detailed theoretical guarantees for our proposed algorithms are provided in Appendix B.

\section{{Notation and Preliminary}} \label{sec:notation}
In this paper, matrices and tensors are represented by capital letters (e.g. $A, B,\dots$) and curly letters (e.g. $\mathcal{A},\mathcal{B},\dots$), respectively. The Matlab command $A'$ can be used to represent the conjugate transpose of the matrix $A$. $\mathbb{R}$ and $\mathbb{C}$ represent the real number space and the complex number space, respectively.
For matrix $A\in \mathbb{C}^{n_1\times n_2}$, its $(i,j)$-th element is represented by $a_{i,j}$. For the third-order tensor $\mathcal{A}\in \mathbb{C}^{n_1\times n_2 \times n_3}$, its $(i,j,k)$-th element is represented by $a_{i,j,k}$. 
The Matlab notations $\mathcal{A}(i,:,:)$, $\mathcal{A}(:,i,:)$ and $\mathcal{A}(:,:, i)$ are used to represent the $i$-th horizontal, lateral and frontal slices of $\mathcal{A}$, respectively. 
The facial slice $\mathcal{A}(:,:,i)$ is also represented by $\mathcal{A}^{(i)}$. The Frobenius norm of a tensor $\mathcal{A}$ is defined as the square root of the sum of the squares of its elements, i.e.,
\begin{equation}
	{\|\mathcal{A}\|} _{F}:={\|\mathcal{A}(:)\|} _{2}=\sqrt{\langle{\mathcal{A},\mathcal{A}}\rangle}=\sqrt{\sum_{ijk}|a_{ijk}|^{2}} \ .
\end{equation}
$\mathcal{A}^{H}$ and $\mathcal{A}^{\dagger}$ represent the conjugate transpose and pseudo-inverse of $\mathcal{A}$, respectively.

This paper focuses on the low tubal rank tensor approximation that meets the desired accuracy in an efficient manner.
For tensor $\mathcal{A} \in \mathbb{R}^{n_{1}\times n_{2}\times n_{3}}$, the mathematical model for finding the low-rank approximation $\hat{\mathcal{A}}$ of $\mathcal{A}$ can be expressed as
\begin{equation}\label{eq2.1}
	{\|\mathcal{A}-\hat{\mathcal{A}}\|} _{F} = \min_{{\rm rank}_{L}(\mathcal{B})\leq l}{\|\mathcal{A}-\mathcal{B}\|} _{F},
\end{equation}
where $l\ll \min\{n_{1},n_{2}\}$ is the target rank, $L$ represents an arbitrary invertible linear transform, and ${\rm rank}_{L}(\mathcal{B})$ denotes the transformed tubal rank of tensor $\mathcal{B}$.

\subsection{{Transformed Tensor SVD}}
We below briefly recall the transformed tensor SVD of third-order tensors; more details can be found in \cite{Kernfeld15}.

For any third-order tensor $\mathcal{A}\in \mathbb{R}^{n_{1}\times n_{2}\times n_{3}}$, let $\bar{\mathcal{A}}_{L}$ represent a third-order tensor obtained via being multiplied by $L$ (an arbitrary invertible linear transform) on all tubes along the third-dimension of $\mathcal{A}$, i.e., 
\begin{equation}
	\bar{\mathcal{A}}_{L}(i,j,:)= L(\mathcal{A}(i,j,:)), \ i=1,\dots,n_{1}, \ j=1,\dots,n_{2}.
\end{equation}
Here we write $\bar{\mathcal{A}}_{L}=L[\mathcal{A}]$. Moreover, one can get $\mathcal{A}$ from $\bar{\mathcal{A}}_{L}$ by using $L^{-1}$ along the third-dimension of $\bar{\mathcal{A}}_{L}$, i.e., $\mathcal{A}=L^{-1}[\bar{\mathcal{A}}_{L}]$.  The $\star{_{L}}$-product is defined in Definition \ref{df2.1} below.

\begin{definition}\label{df2.1}
	\cite{Kernfeld15} For any two tensors $\mathcal{X}\in \mathbb{C}^{n_{1}\times n_{2}\times n_{3}}$ and $\mathcal{Y}\in \mathbb{C}^{n_{2}\times n_{4}\times n_{3}}$, and an arbitrary invertible linear transform $L$, the $\star{_{L}}$-product of $\mathcal{X}$ and $\mathcal{Y}$ is a tensor $\mathcal{Z}\in \mathbb{C}^{n_{1}\times n_{4}\times n_{3}}$ given by
	\begin{equation}		\mathcal{Z}=\mathcal{X}\star{_{L}}\mathcal{Y}=L^{H}[\text{fold}(\text{block}(\bar{\mathcal{X}}_{L})\text{block}(\bar{\mathcal{Y}}_{L}))],
	\end{equation}
	where $\text{fold}(\text{block}(\bar{\mathcal{X}}_{L}))=\bar{\mathcal{X}}_{L}$ and $\bar{X}_{L}=\text{block}(\bar{\mathcal{X}}_{L})=\begin{pmatrix} \bar{\mathcal{X}}_{L}^{(1)} & & &  \\  & \bar{\mathcal{X}}_{L}^{(2)}& & \\ & &\ddots&\\  & & &\bar{\mathcal{X}}_{L}^{(n_{3})}\end{pmatrix}$. 
\end{definition}

\begin{definition}\label{df2.2}
	The Kronecker product of matrices $A\in \mathbb{C}^{m\times l}$ and $B\in \mathbb{C}^{p\times r}$ is
	$$A \otimes B=\begin{pmatrix} a_{11}B &a_{12}B  &\dots&a_{1l}B \\ a_{21}B & a_{22}B &\dots&a_{2l}B \\ \vdots &\vdots &\ddots&\vdots\\ a_{m1}B & a_{m2}B&\dots&a_{ml}B\end{pmatrix}. $$
\end{definition}

The t-product \cite{Kilmer11} of $\mathcal{A}\in \mathbb{C}^{n_{1}\times n_{2}\times n_{3}}$ and $\mathcal{B}\in \mathbb{C}^{n_{2}\times n_{4}\times n_{3}}$ is a tensor $\mathcal{C}\in \mathbb{C}^{n_{1}\times n_{4}\times n_{3}}$ given by
\begin{equation}\label{eq2.2}
	\begin{split}	\mathcal{C}=\mathcal{A}*\mathcal{B}=\text{fold}_{\text{vec}}(\text{circ}(\mathcal{A})\times \text{vec}(\mathcal{B})),
	\end{split}
\end{equation}
where
$$\text{vec}(\mathcal{B})=\begin{pmatrix} \mathcal{B}^{(1)} \\\mathcal{B}^{(2)} \\  \vdots\\\mathcal{B}^{(n_{3})}\end{pmatrix} \in \mathbb{R}^{n_{2}n_{3}\times n_{4}}, \quad \text{fold}_{\text{vec}}(\text{vec}(\mathcal{B}))=\mathcal{B}, $$
$$\text{circ}(\mathcal{A})=\begin{pmatrix} \mathcal{A}^{(1)} &\mathcal{A}^{(n_{3})}  &\mathcal{A}^{(n_{3}-1)} &\dots&\mathcal{A}^{(2)} \\ \mathcal{A}^{(2)} & \mathcal{A}^{(1)} &\mathcal{A}^{(n_{3})}&\dots&\mathcal{A}^{(3)} \\ \vdots &\vdots &\vdots &\ddots&\vdots\\ \mathcal{A}^{(n_{3})} & \mathcal{A}^{(n_{3}-1)}&\mathcal{A}^{(n_{3}-2)} &\dots&\mathcal{A}^{(1)}\end{pmatrix}.$$

Let $\text{$F_{n_{3}}$}$ and $C$ represent the discrete Fourier transform (DFT) matrix and the discrete cosine transform (DCT) matrix, respectively. The t-product in Eq. \eqref{eq2.2} can be seen as a special case of Definition \ref{df2.1}. Recall that the block circulant matrix $\text{circ}(\mathcal{A})$ can be
diagonalized by the fast Fourier transform matrix $\text{$F_{n_{3}}$}$, and the block diagonal matrices are the frontal slices of $\bar{\mathcal{A}}_\text{$F_{n_{3}}$}$,
i.e.,
\begin{equation}
	\bar{A}_\text{$F_{n_{3}}$} 
	= \text{block}(\bar{\mathcal{A}}_\text{$F_{n_{3}}$})
	= (\text{$F_{n_{3}}$}\otimes I_{n_{1}})\times \text{circ}(\mathcal{A})\times (\text{$F^{H}_{n_{3}}$}\otimes I_{n_{2}}).
\end{equation}
It follows that
\begin{eqnarray*}
	\mathcal{A}*\mathcal{B} & = & \text{fold}_{\text{vec}}(\text{circ}(\mathcal{A})\times \text{vec}(\mathcal{B})) \\
	& = &  \text{fold}_{\text{vec}}((\text{$F^{H}_{n_{3}}$}\otimes I_{n_{1}})\times \text{block}(\bar{\mathcal{A}}_\text{$F_{n_{3}}$}) \times (\text{$F_{n_{3}}$}\otimes I_{n_{2}})\times \text{vec}(\mathcal{B}))\\
	& = & \text{fold}((\text{$F^{H}_{n_{3}}$}\otimes I_{n_{1}})\times \text{block}(\bar{\mathcal{A}}_\text{$F_{n_{3}}$}) \times \text{block}(\bar{\mathcal{B}}_\text{$F_{n_{3}}$}))\\
	& = & \text{$F^{H}_{n_{3}}$}[\text{fold}(\text{block}(\bar{\mathcal{A}}_\text{$F_{n_{3}}$})\text{block}(\bar{\mathcal{B}}_\text{$F_{n_{3}}$}))]\\
	& = & \text{$F^{H}_{n_{3}}$}[\text{fold}(\bar{A}_\text{$F_{n_{3}}$}\bar{B}_\text{$F_{n_{3}}$})]\\
	& = & \mathcal{A}\star_{\text{$F_{n_{3}}$}}\mathcal{B} .
\end{eqnarray*}

The definitions of the conjugate transpose of tensor, the identity tensor, the unitary tensor, the invertible tensor, and the diagonal tensor related to the $\star{_{L}}$-product are given as follows.
\begin{itemize}
	\item \cite{Song20} The conjugate transpose of $\mathcal{A}\in \mathbb{C}^{n_{1}\times n_{2}\times n_{3}}$ with respect to $L$ is the tensor $\mathcal{A}^{H}\in \mathbb{C}^{n_{1}\times n_{2}\times n_{3}}$ obtained by $\mathcal{A}^{H}=L^{H}[\text{fold}(\text{block}(\bar{\mathcal{A}}_{L})^{H})]=L^{H}[\text{fold}(\bar{A}^{H}_{L})].$
	\item \cite{Kernfeld15} The identity tensor $\mathcal{I}_{L}\in \mathbb{C}^{n\times n\times n_{3}}$ (with respect to $L$) is defined to be a
	tensor such that $\mathcal{I}_{L}=L^{H}[\mathcal{I}]$, where each frontal slice of $\mathcal{I}$ is the $n\times n$ identity matrix.
	\item \cite{Kernfeld15} A tensor $\mathcal{Q}\in \mathbb{C}^{n\times n\times n_{3}}$ is unitary with respect to $\star{_{L}}$-product if it satisfies ${\mathcal{Q}^{H}}\star{_{L}}\mathcal{Q}=\mathcal{Q}
	\star{_{L}}\mathcal{Q}^{H}=\mathcal{I}_{L}$, where $\mathcal{I}_{L}$ is the identity tensor.
	\item \cite{Kernfeld15} For tensors $\mathcal{A}\in \mathbb{C}^{n\times n\times n_{3}}$ and $\mathcal{B}\in \mathbb{C}^{n\times n\times n_{3}}$, if ${\mathcal{A}}\star{_{L}}\mathcal{B}=\mathcal{B}\star{_{L}}\mathcal{A}=\mathcal{I}_{L}$, then tensor $\mathcal{B}$ is the invertible tensor under the $\star{_{L}}$-product of tensor $\mathcal{A}$.
	\item \cite{Kilmer11} A tensor is a diagonal tensor if each frontal slice of the tensor is a diagonal matrix. For a third-order tensor, if all of its frontal slices are upper or lower triangles, then the tensor is called f-upper or f-lower.
\end{itemize}

\begin{lemma}\label{lem2.3}
	\cite{Kernfeld15, Kilmer11} Suppose that $\mathcal{A}\in \mathbb{R}^{m\times k \times p}$ and $\mathcal{B}\in \mathbb{R}^{k\times n \times p}$ are two arbitrary tensors. Let $\mathcal{Z}=\mathcal{A}\star{_{L}}\mathcal{B}$. Then the following properties hold:
	\begin{itemize}
		\item[(1)] ${{\|\mathcal{A}\|} _{F}^{2}}=\frac{1}{p}{\|\bar{A}_{\text{$F_{n_{3}}$}}\|} _{F}^{2}=\frac{1}{p}{\|\bar{A}_{L}\|} _{F}^{2}=\frac{1}{p}\sum_{i=1}^{p}{
			\|\bar{\mathcal{A}}^{(i)}_{L}\|} _{F}^{2}$.
		\item[(2)] $\mathcal{Z}=\mathcal{A}\star{_{L}}\mathcal{B}$ is equivalent to $\text{block}(\bar{\mathcal{Z}}_{L})=\text{block}(\bar{\mathcal{A}}_{L})\text{block}(\bar{\mathcal{B}}_{L})$$, i.e., \bar{Z}_{L}=\bar{A}_{L}\bar{B}_{L}$.
	\end{itemize}
\end{lemma}

\begin{definition}\label{df2.4}
	(Gaussian Random Tensor) \cite{Zhang18} A tensor $\mathcal{G}\in \mathbb{R}^{m\times n\times p}$
	is called a Gaussian random tensor if the elements of $\mathcal{G}^{(1)}$ satisfy the standard normal distribution
	(i.e., Gaussian with mean zero and variance one) and the other frontal slices are all zero.
\end{definition}

Based on the above
definitions, we have the following tensor SVD with respect to $L$.
\begin{theorem}\cite{Kernfeld15}\label{thm2.5}
	$\forall \mathcal{A}\in \mathbb{C}^{m\times n\times p}$, the transformed tensor SVD is given by
	\begin{equation}\label{eq2.3}	\mathcal{A}=\mathcal{U}\star{_{L}}\mathcal{S}\star{_{L}}\mathcal{V}^{H},
	\end{equation}
	where $\mathcal{U}\in \mathbb{C}^{m\times m\times p}$ and $\mathcal{V}\in \mathbb{C}^{n\times n\times p}$ are unitary tensors with respect to the $\star{_{L}}$-product, and $\mathcal{S}\in \mathbb{C}^{m\times n\times p}$ is a diagonal tensor.
\end{theorem}

The transformed tubal rank, denoted as
${\rm rank}_{L}(\mathcal{A})$, is defined as the number of nonzero singular tubes of $\mathcal{S}$, 
i.e.,
\begin{equation}\label{eq2.4}
	\begin{split}
		{\rm rank}_{L}(\mathcal{A})=\# \{i:\mathcal{S}(i,i,:)\ne 0\},
	\end{split}
\end{equation}
where \# denotes the cardinality of a set.
The transformed tensor SVD could be implemented efficiently by the SVDs of the frontal slices
in the transformed domain. We also refer the readers to \cite{Song20} for more details on the
computation of the transformed tensor SVD.
For the Kernfeld-Kilmer transformed tensor SVD (i.e., Eq. \eqref{eq2.3}), by \cite{Kernfeld15, Kilmer11, Kilmer08}, we have
\begin{equation}\label{eq2.5}
	\sum_{k=1}^{p}\mathcal{S}(1,1,k)^{2}\geq \sum_{k=1}^{p}\mathcal{S}(2,2,k)^{2}\geq \dots \geq \sum_{k=1}^{p}\mathcal{S}(\min\{m,n\}, \min\{m,n\},k)^{2}.
\end{equation}

\begin{definition}\label{df2.6}
	(Transformed Tensor Singular Values) 
	Suppose $\mathcal{A}\in \mathbb{R}^{m\times n\times p}$ with a Kernfeld-Kilmer transformed tensor SVD such that Eq. \eqref{eq2.5} is satisfied. The $i$-th largest transformed tensor singular value of $\mathcal{A}$ is defined as
	\begin{equation}\label{eq2.6}
		\sigma_{i}=\sqrt{\sum_{k=1}^{p}\mathcal{S}(i,i,k)^{2}}, \quad
		{\rm for} \ i=1,2,\dots, \min\{m,n\}.
	\end{equation}
	
\end{definition}

Similarly to the definition of the matrix tail energy, the tail energy of tensor is defined below.
\begin{definition}\label{df2.7}
	(Tail Energy) For tensor $\mathcal{A}\in \mathbb{R}^{m\times n\times p}$, the $i$-th largest transformed tensor singular value is $\sigma_{i}$, $i=1,\dots, \min\{m,n\}$. Then the $j$-th tail energy of $\mathcal{A}$ is defined as
	\begin{equation}\label{eq2.7}
		{\tau_{j}}^{2}(\mathcal{A}):=\min_{{\rm rank}_{L}(\mathcal{B})<j}{\|\mathcal{A}-\mathcal{B}\|}_{F}^{2}=\sum_{i\ge j}\sigma_{i}^{2}(\mathcal{A}) .
	\end{equation}	
\end{definition}

According to the above definition and using Lemma \ref{lem2.3} and the linearity, we can obtain
the following proposition.
\begin{proposition} \label{pro2.8}
	Suppose $\mathcal{A}\in \mathbb{R}^{m\times n\times p}$ and $\bar{\mathcal{A}}_{L}$ represents a third-order tensor obtained via being multiplied by $L$ on all tubes along the third-dimension of $\mathcal{A}$. Let $j$ be a positive integer satisfying $j\le \min\{m,n\}$. Then
	\begin{equation}\label{eq2.8}
		{\tau_{j}}^{2}(\mathcal{A})=\frac{1}{p}\sum_{i=1}^{p}{\tau_{j}}^{2}(\bar{\mathcal{A}}_{L}^{(i)}) .
	\end{equation}	
\end{proposition}


\subsection{\textbf{Tensor Sketching Operator}} \label{subsec:tso}
Using the three matrix sketching techniques (i.e., Gaussian projection, subsampled randomized Hadamard transform (SRHT), and count sketch \cite{{Wang15}}) with their pseudocodes (i.e., \texttt{GaussianProjection}, \texttt{SRHT}, and \texttt{CountSketch}) introduced in Appendix A, three corresponding tensor sketching operators can be generated. As for an efficient two-sided sketching algorithm, we need to generate four random linear dimension reduction maps, i.e.,
\begin{equation}\label{eq2.9}
	\displaystyle \Upsilon \in \mathbb{R}^{k\times m\times p}, \ \  \Omega \in \mathbb{R}^{k\times n \times p}, \ \ 
	\displaystyle \Phi \in \mathbb{R}^{s\times m\times p}, \ \  {\rm and} \ \Psi \in \mathbb{R}^{s\times n \times p} .
\end{equation}
Different ways of generating tensor sketch operators regarding $\mathcal{A}\in \mathbb{R}^{m\times n\times p}$ are shown below.
\begin{itemize}
	\item \textbf{Gaussian tensor sketching operator.}
	Set\\
	$\Upsilon=\text{zeros}(k,m,p);  \ \Upsilon(:,:,1)=\texttt{GaussianProjection}(\mathcal{A}(:,:,1),k)'$;\\
	$\Omega=\text{zeros}(k,n,p);\ \Omega(:,:,1)=\texttt{GaussianProjection}(\mathcal{A}(:,:,1)',k)';\\
	\Phi=\text{zeros}(s,m,p); \ \Phi(:,:,1)=\texttt{GaussianProjection}(\mathcal{A}(:,:,1),s)'$;\\
	$\Psi=\text{zeros}(s,n,p); \ \Psi(:,:,1)=\texttt{GaussianProjection}(\mathcal{A}(:,:,1)',s)'$.\\ 
	Then $\Upsilon$, $\Omega$, $\Phi$ and $\Psi$ are said to be the Gaussian tensor sketching operators. 
	
	\item 	\textbf{SRHT tensor sketching operator.}
	Set\\ 
	$\Upsilon=\text{zeros}(k,m,p);\ \Upsilon(:,:,1)=\texttt{SRHT}(\mathcal{A}(:,:,1),k)'$;\\
	$\Omega=\text{zeros}(k,n,p);$ \ 
	$\Omega(:,:,1)=\texttt{SRHT}(\mathcal{A}(:,:,1)',k)';\\
	\Phi=\text{zeros}(s,m,p);\ \Phi(:,:,1)=\texttt{SRHT}(\mathcal{A}(:,:,1),s)'$;\\
	$\Psi=\text{zeros}(s,n,p);\ \Psi(:,:,1)=\texttt{SRHT}(\mathcal{A}(:,:,1)',s)'$. \\
	Then $\Upsilon$, $\Omega$, $\Phi$ and $\Psi$ are said to be the SRHT tensor sketching operators.
	
	\item \textbf{ Count tensor sketching operator.}
	For $i = 1,2,\ldots, p$, set\\
	$\Upsilon=\text{zeros}(k,m,p); \ \Upsilon(:,:,i)=\texttt{CountSketch}(\mathcal{A}(:,:,1), k)'$;\\
	$\Omega=\text{zeros}(k,n,p); \ \Omega(:,:,i)=\texttt{CountSketch}(\mathcal{A}(:,:,1)', k)';\\
	\Phi=\text{zeros}(s,m,p); \ \Phi(:,:,i)=\texttt{CountSketch}(\mathcal{A}(:,:,1), s)'$;\\
	$\Psi=\text{zeros}(s,n,p); \ \Psi(:,:,i)=\texttt{CountSketch}(\mathcal{A}(:,:,1)', s)'$.\\
	Then $\Upsilon$, $\Omega$, $\Phi$ and $\Psi$ are said to be the count tensor sketching operators.
\end{itemize}

\section{The Proposed Two-Sided Sketching Algorithms} \label{sec:our-ts-alg}
The two-sided tensor sketching algorithm proposed by Qi and Yu \cite{Qi21} only considers the range and co-range of the input tensor. On this basis, we here also consider the core sketch. The core sketch contains new information that improves our estimates of the transformed tensor singular values and the transformed tensor singular vectors of the input tensor, and is responsible for the superior performance of the algorithms. We below firstly present the framework of our efficient two-sided sketching method based on the transformed domains for low tubal rank tensor approximation, followed by the principle interpreting its rationale and its extension by using the power iteration technique.

\subsection{Method}
Given the input tensor $\mathcal{A}\in \mathbb{R}^{m\times n\times p}$ and the objective tubal rank $k$, using the appropriate tensor sketching operators $\Upsilon, \Omega, \Phi$ and $\Psi$ in Eq. \eqref{eq2.9}, we can realize the randomized sketches $(\mathcal{X},\mathcal{Y},\mathcal{Z})$  such as
\begin{align}
	\mathcal{X}& :=\Upsilon \star{_{L}} \mathcal{A}\in \mathbb{R}^{k\times n \times p}, \label{eqn:rs-x} \\
	\mathcal{Y} & :=\mathcal{A} \star{_{L}} {\Omega}^{H}\in \mathbb{R}^{m\times k \times p}, \label{eqn:rs-y} \\
	\mathcal{Z}& :=\Phi \star{_{L}} \mathcal{A} \star{_{L}} \Psi^{H} \in  \mathbb{R}^{s\times s \times p}. 
	\label{eqn:rs-z}
\end{align}
The first two tensor sketches $\mathcal{X}$ and $\mathcal{Y}$ respectively capture the co-range and range of $\mathcal{A}$, and the core sketch $\mathcal{Z}$, as we mentioned before, contains new information that improves our estimates of the transformed tensor singular values and vectors of $\mathcal{A}$ and is also responsible for further method performance enhancement.

Once the sketches $(\mathcal{X},\mathcal{Y},\mathcal{Z})$ of the input tensor $\mathcal{A}$ are obtained by Eq. \eqref{eqn:rs-x}-- \eqref{eqn:rs-z}, we can find the low-rank approximation $\hat{\mathcal{A}}$ by following the below proposed three-step process.

\begin{itemize}
	\item[(I)] Form an L-orthogonal-triangular factorization
	\begin{equation}\label{eq3.12}
		\mathcal{X}^{H}:=\mathcal{P} \star{_{L}} \mathcal{R}_{1}, \quad 
		\mathcal{Y}:=\mathcal{Q} \star{_{L}} \mathcal{R}_{2},
	\end{equation}
	where $\mathcal{P} \in \mathbb{R}^{n\times k \times p}$ and $\mathcal{Q} \in \mathbb{R}^{m\times k \times p}$
	are  partially orthogonal tensors, and $\mathcal{R}_{1}\in \mathbb{R}^{k\times k \times p}$ and $\mathcal{R}_{2}\in \mathbb{R}^{k\times k \times p}$
	are f-upper
	triangular tensors, in the sense of the $\star{_{L}}$-product operation.
	
	\item[(II)]  
	Using the known $\mathcal{P}, \mathcal{Q}, \Phi, \Psi$, and $\mathcal{Z}$, calculate 
	\begin{equation}\label{eq3.15}
		\mathcal{C}:=(\Phi \star{_{L}} \mathcal{Q})^{\dagger}\star{_{L}} \mathcal{Z}\star{_{L}} ((\Psi \star{_{L}} \mathcal{P})^{\dagger})^{H}\in \mathbb{R}^{k\times k \times p} \ .
	\end{equation}
	The above formula of calculating $\mathcal{C}$ is equivalent to solving the below least-squares problem based on the $\star{_{L}}$-product, i.e.,
	\begin{equation} \label{eqn:ls-4-C}
		\min_{\cal C} \frac{1}{2} {\| \mathcal{G}\star{_{L}} \mathcal{C}\star{_{L}} \mathcal{M}^{H}-\mathcal{Z} \|}_{F}^{2} \ ,
	\end{equation}
	where $\mathcal{G} = \Phi \star{_{L}} \mathcal{Q}$ and $\mathcal{M} = \Psi \star{_{L}} \mathcal{P}$.  According to Lemma \ref{lem2.3}, the above problem \eqref{eqn:ls-4-C} can be reformulated as
	\begin{equation}\label{eq3.14}
		\min_{\bar{\mathcal{C}}_{L}} \frac{1}{2p} {\| \bar{\mathcal{G}}_{L} \bar{\mathcal{C}}_{L} \bar{\mathcal{M}}^{H}_{L}-\bar{\mathcal{Z}}_{L} \|}_{F}^{2} \ ,
	\end{equation}
	whose solution is $\bar{\mathcal{C}}^{(i)}_{L}=({\bar{\mathcal{G}}^{(i)}_{L}})^{\dagger} \bar{\mathcal{Z}}^{(i)}_{L}((\bar{\mathcal{M}} ^{(i)}_{L})^{\dagger})^{H}$ for $i = 1, \ldots, p$. 
	
	\item[(III)] 
	Construct the transformed tensor tubal rank $k$ approximation
	\begin{equation}\label{eqn:app-A-final}
		\hat{\mathcal{A}}:=\mathcal{Q}\star{_{L}} \mathcal{C}\star{_{L}} {\mathcal{P}}^{H} \ .
	\end{equation}
\end{itemize}

\begin{algorithm}[htb] 
	\caption{$L$ Transformed Randomized Projection Sketching Algorithm}
	\begin{algorithmic}[1] 
		\STATE {\textbf{Input:}} Input tensor $\mathcal{A}\in \mathbb{R}^{m\times n\times p}$ and sketch size parameters $k, s$ with $k\leq s$.
		\STATE {\textbf{function}} \texttt{L-TRP-SKETCH}{($\mathcal{A},k$)}
		\STATE Select the appropriate tensor sketching operators from Eq. (\ref{eq2.9}), i.e., \\ $\Upsilon \in \mathbb{R}^{k\times m\times p}, \Omega \in \mathbb{R}^{k\times n \times p}, \Phi \in \mathbb{R}^{s\times m\times p}, \Psi \in \mathbb{R}^{s\times n \times p}$;
		\STATE ${\bar{\mathcal{A}}}_{L} = L[\mathcal{A}]$, ${\bar{\Upsilon}}_{L} = L[\Upsilon]$, ${\bar{\Omega}}_{L} = L[\Omega]$, ${\bar{\Phi}}_{L} = L[\Phi]$, ${\bar{\Psi}}_{L} = L[\Psi]$;
		\STATE {\textbf{for}} {$i \gets 1 \;\textbf{to} \;p$}
		\STATE \quad ${\bar{\mathcal{X}}^{(i)}}_{L} = {\bar{\Upsilon}^{(i)}}_{L}{\bar{\mathcal{A}}^{(i)}}_{L}$, 
		${\bar{\mathcal{Y}}}^{(i)}_{L} = {\bar{\mathcal{A}}}^{(i)}_{L}({\bar{\Omega}}^{(i)}_{L})^{H}$,  
		${\bar{\mathcal{Z}}}^{(i)}_{L} = {\bar{\Phi}}^{(i)}_{L}({\bar{\mathcal{A}}}^{(i)}_{L})({\bar{\Psi}}^{(i)}_{L})^{H}$;
		\STATE  \quad $[\bar{\mathcal{P}}^{(i)}_{L},\bar{\mathcal{R}}_{L}^{(i)}]=\texttt{qr}(({\bar{\mathcal{X}}}^{(i)}_{L})^{H},0),  \ [\bar{\mathcal{Q}}^{(i)}_{L},\bar{\mathcal{M}}_{L}^{(i)}]=\texttt{qr}({\bar{\mathcal{Y}}}^{(i)}_{L},0)$;
		\STATE \quad $\bar{\mathcal{C}}^{(i)}_{L}=(\bar{\Phi}^{(i)}_{L} \bar{\mathcal{Q}}^{(i)}_{L})^{\dagger} \bar{\mathcal{Z}}^{(i)}_{L}((\bar{\Psi} ^{(i)}_{L} \bar{\mathcal{P}}^{(i)}_{L})^{\dagger})^{H}$;
		\STATE  \quad $\tilde{\mathcal{A}}^{(i)}_{L}=\bar{\mathcal{Q}}^{(i)}_{L} \bar{\mathcal{C}}^{(i)}_{L} {(\bar{\mathcal{P}}^{(i)}_{L})}^{H}$;
		\STATE {\textbf{end}}
		\STATE {\textbf{return}} \ $\hat{\mathcal{A}}_{L} =L^{H}[{\tilde{\mathcal{A}}}_{L}]$ and  $({\bar{\mathcal{Q}}}_{L}, {\bar{\mathcal{C}}}_{L}, {\bar{\mathcal{P}}}_{L})$.
		
	\end{algorithmic}
	\label{alg:l-trp-sketch}
\end{algorithm}


We call the above three-step process our proposed two-sided sketching method, i.e., {\it $L$ transformed randomized projection sketching (L-TRP-SKETCH)} algorithm. The pseudocode of our L-TRP-SKETCH algorithm is given in Algorithm \ref{alg:l-trp-sketch}. The storage cost for the sketches $(\mathcal{X},\mathcal{Y},\mathcal{Z})$ is $p(nk+mk+s^{2})$ floating point numbers. The storage complexity of  Algorithm 3.1
for the original data $\mathcal{A}\in \mathbb{R}^{m\times n\times p}$ is $\mathcal{O}(mkp+nkp+kkp) $ floating point numbers.
Regarding the selection of the invertible transformation operator $L$, in addition to the DFT matrix $F_{n_3}$ (here $n_3 = p$) and the DCT matrix $C$, we can also choose the U transformation matrix as in \cite{Song20}. In this case, the original data $\mathcal{A}\in \mathbb{R}^{m\times n\times p}$ modulo-3 is expanded into a matrix $W\in \mathbb{R}^{p\times mn}$, and then SVD is applied on $W$ to obtain the unitary transformation matrix $U\in \mathbb{R}^{p\times p}$.

\textit{Naming.} When, based on the $U$ transformed domain, linear dimensionality reduction mappings are chosen as the count tensor sketching operator, the Gaussian tensor sketching operator, and the SRHT tensor sketching operator, Algorithm \ref{alg:l-trp-sketch} is then referred to as the U-Count-Sketch, the U-Gaussian-Sketch, and the U-SRHT-Sketch algorithms, respectively. 
When, based on the DFT transformed domain, linear dimensionality reduction mappings are chosen as the count tensor sketching operator, the Gaussian tensor sketching operator, and the SRHT tensor sketching operator, Algorithm \ref{alg:l-trp-sketch} is then referred to as the DFT-Count-Sketch, the DFT-Gaussian-Sketch, and the DFT-SRHT-Sketch algorithms, respectively.
Similarly, when, based on the DCT transformed domain, linear dimensionality reduction mappings are chosen as the count tensor sketching operator, the Gaussian tensor sketching operator, and the SRHT tensor sketching operator, Algorithm \ref{alg:l-trp-sketch} is then referred to as the DCT-Count-Sketch, the DCT-Gaussian-Sketch, and the DCT-SRHT-Sketch algorithms, respectively. Finally, if the transformed domain is not specified, we just use L in the names of these algorithms; e.g., U-Gaussian-Sketch will be called L-Gaussian-Sketch.

\subsection{Principle}
For 
\begin{equation} \label{eq3.17}
	\mathcal{A}\approx \mathcal{Q}\star_{L}(\mathcal{Q}^{H}\star_{L}\mathcal{A}\star_{L}\mathcal{P})\star_{L} \mathcal{P}^{H},
\end{equation}
the core tensor $\mathcal{Q}^{H}\star_{L}\mathcal{A}\star_{L}\mathcal{P}$ cannot be calculated directly from the linear sketch since $\mathcal{P}$ and $\mathcal{Q}$ are functions of $\mathcal{A}$. Using the representation in Eq. (\ref{eq3.17}), the core sketch $\mathcal{Z}$ estimating the core tensor can be achieved by
\begin{equation} \label{eq3.18}
	\mathcal{Z} = \Phi \star_{L} \mathcal{A}\star_{L}{\Psi}^{H} \approx (\Phi \star_{L}\mathcal{Q})\star_{L}(\mathcal{Q}^{H}\star_{L}\mathcal{A}\star_{L}
	\mathcal{P})\star_{L}(\mathcal{P}^{H}\star_{L}\Psi^{H}).
\end{equation}
Transferring the external matrix to the left-hand side, the core approximation $\mathcal{C}$ defined in Eq. (\ref{eq3.15}) is found to satisfy
\begin{equation}\label{eq3.19}
	\mathcal{C}=(\Phi \star_{L}\mathcal{Q})^{\dagger}\star_{L}\mathcal{Z}\star_{L} ((\Psi \star_{L}\mathcal{P})^{\dagger})^{H}\approx \mathcal{Q}^{H}\star_{L}\mathcal{A}\star_{L}\mathcal{P}.
\end{equation}
Given Eq. \eqref{eqn:app-A-final}, (\ref{eq3.17}) and (\ref{eq3.19}), we have
\begin{equation} \label{eq3.20}
	\mathcal{A}\approx \mathcal{Q}\star_{L}(\mathcal{Q}^{H}\star_{L}\mathcal{A}\star_{L}\mathcal{P})\star_{L}
	\mathcal{P}^{H}\approx \mathcal{Q}\star_{L}\mathcal{C}\star_{L}\mathcal{P}^{H} = \hat{\mathcal{A}}.
\end{equation}
The error in the last relation depends on the error in the best transformed tubal rank $k$ approximation of $\mathcal{A}$.

\begin{algorithm}[htb] 
	\caption{DCT-Gaussian-Sketch-PI}
	\begin{algorithmic}[1] 
		\STATE {\textbf{Input:}} Input tensor $\mathcal{A}\in \mathbb{R}^{m\times n\times p}$ and sketch size parameters $k, s, q$ with $k\leq s$.
		
		\STATE {\textbf{function}} \texttt{DCT-Gaussian-Sketch-PI}{($\mathcal{A},k,q$)}
		\STATE Select the Gaussian tensor sketching operators from Eq. (\ref{eq2.9}), i.e., \\
		$\Upsilon \in \mathbb{R}^{k\times m\times p}, \Omega \in \mathbb{R}^{k\times n \times p}, \Phi \in \mathbb{R}^{s\times m\times p}, \Psi \in \mathbb{R}^{s\times n \times p}$;
		\STATE ${\bar{\mathcal{A}}}_{C} = C[\mathcal{A}]$; ${\bar{\Upsilon}}_{C} = C[\Upsilon]$; ${\bar{\Omega}}_{C} = C[\Omega]$; ${\bar{\Phi}}_{C} = C[\Phi]$; ${\bar{\Psi}}_{C} = C[\Psi]$;
		
		\STATE {\textbf{for}} {$i \gets 1 \;\textbf{to} \;p$}
		
		\STATE \quad 
		${\bar{\mathcal{X}}^{(i)}}_{C} = {\bar{\Upsilon}^{(i)}}_{C}{\bar{\mathcal{A}}^{(i)}}_{C}$;
		${\bar{\mathcal{Y}}}^{(i)}_{C} = {\bar{\mathcal{A}}}^{(i)}_{C}({\bar{\Omega}}^{(i)}_{C})^{H}$;
		${\bar{\mathcal{Z}}}^{(i)}_{C} = {\bar{\Phi}}^{(i)}_{C}({\bar{\mathcal{A}}}^{(i)}_{C})({\bar{\Psi}}^{(i)}_{C})^{H}$;
		
		\STATE  \quad 
		$[\bar{\mathcal{P}}^{(i)}_{C}, \bar{\mathcal{R}}_{C}^{(i)}]=\texttt{qr}(({\bar{\mathcal{X}}}^{(i)}_{C})^{H},0), [\bar{\mathcal{Q}}^{(i)}_{C}, \bar{\mathcal{M}}_{C}^{(i)}]=\texttt{qr}({\bar{\mathcal{Y}}}^{(i)}_{C},0)$;

		\STATE \quad {\textbf{for}} {$j \gets 1 \;\textbf{to} \;q$}

		\STATE \quad \quad${\tilde{\mathcal{Y}}}^{(i)}_{C}=({\bar{\mathcal{A}}^{(i)}}_{C})^{H}{\bar{\mathcal{Q}}^{(i)}}_{C},
		[\tilde{\mathcal{Q}}^{(i)}_{C}, \sim]=\texttt{qr}({\tilde{\mathcal{Y}}}^{(i)}_{C},0)$;
		\STATE \quad \quad
		${\hat{\mathcal{Y}}}^{(i)}_{C}={\bar{\mathcal{A}}^{(i)}}_{C}{\tilde{\mathcal{Q}}^{(i)}}_{C},
		[\hat{\mathcal{Q}}^{(i)}_{C}, \sim]=\texttt{qr}({\hat{\mathcal{Y}}}^{(i)}_{C},0)$;
		\STATE \quad \quad
		${\tilde{\mathcal{X}}^{(i)}}_{C} = {\bar{\mathcal{A}}^{(i)}}_{C} {\bar{\mathcal{P}}^{(i)}}_{C},
		[\tilde{\mathcal{P}}^{(i)}_{C}, \sim]=\texttt{qr}({\tilde{\mathcal{X}}}^{(i)}_{C},0)$;
		
		\STATE \quad \quad
		${\hat{\mathcal{X}}^{(i)}}_{C} = ({{\bar{\mathcal{A}}^{(i)}}_{C}})^{H}{\tilde{\mathcal{P}}^{(i)}}_{C},
		[\hat{\mathcal{P}}^{(i)}_{C}, \sim]=\texttt{qr}({\hat{\mathcal{X}}}^{(i)}_{C},0)$;
		
		\STATE \quad \quad
		$\bar{\mathcal{Q}}^{(i)}_{C}=\hat{\mathcal{Q}}^{(i)}_{C}, \bar{\mathcal{P}}^{(i)}_{C}=\hat{\mathcal{P}}^{(i)}_{C}$;
		
		\STATE \quad  {\textbf{end}}
		\STATE  \quad $\bar{\mathcal{C}}^{(i)}_{C}=(\bar{\Phi}^{(i)}_{C} \bar{\mathcal{Q}}^{(i)}_{C})^{\dagger} \bar{\mathcal{Z}}^{(i)}_{C}((\bar{\Psi} ^{(i)}_{C} \bar{\mathcal{P}}^{(i)}_{C})^{\dagger})^{H}$;
		\STATE   \quad $\tilde{\mathcal{A}}^{(i)}_{C}=\bar{\mathcal{Q}}^{(i)}_{C} \bar{\mathcal{C}}^{(i)}_{C} {(\bar{\mathcal{P}}^{(i)}_{C})}^{H}$;
		\STATE {\textbf{end}}
		\STATE {\textbf{return}} \ 
		$\hat{\mathcal{A}} =C^{H}[{\tilde{\mathcal{A}}}_{C}]$ and $({\bar{\mathcal{Q}}}_{C}, {\bar{\mathcal{C}}}_{C}, {\bar{\mathcal{P}}}_{C})$. 
		
	\end{algorithmic}
	\label{alg:DCT-Sketch-PI}
\end{algorithm}

\subsection{Extension}
We now showcase one extension of the proposed two-sided sketching algorithms by exploiting the power iteration technique.
As shown in \cite{Halko11}, the power iteration technique is useful to improve sketching algorithms for low-rank matrix approximation. Here, as an example, we can combine the power iteration technique with the DCT-Gaussian-Sketch algorithm, in which we exploit the third order tensor say $\mathcal{B}=(\mathcal{A} \star{_{L}} \mathcal{A}^H)^q \star{_{L}} \mathcal{A}$ (where $q$ is a nonnegative integer) instead of the original tensor $\mathcal{A}$; and the DCT-Gaussian-Sketch algorithm is applied to the new tensor $\mathcal{B}$. According to the transformed tensor SVD, i.e., $\mathcal{A} = \mathcal{U} \star{_{L}} \mathcal{S} \star{_{L}} \mathcal{V}^H$, we have 
\begin{equation}
	\mathcal{B}=\mathcal{U} \star{_{L}} (\mathcal{S})^{2q+1} \star{_{L}} \mathcal{V}^H.
\end{equation}
Therefore, the transformed tensor singular values of $\mathcal{B}$ have a faster decay rate. This way can improve the solution obtained by the DCT-Gaussian-Sketch algorithm. The scheme of the proposed algorithm DCT-Gaussian-Sketch with power iteration, named DCT-Gaussian-Sketch-PI, is summarized in Algorithm \ref{alg:DCT-Sketch-PI}.

\section{Theoretical Analysis}
\label{sec:err-ana}
The error bound of the proposed two-sided sketching algorithms is given in Theorems \ref{thm4.7} and \ref{thmPI} below. The detailed proofs used in the proof of Theorems \ref{thm4.7} and \ref{thmPI} can be found in Appendix B.

\begin{theorem}\label{thm4.7}
	Assume that the sketch parameters satisfy $s\ge 2k+1$ and $\hat{\mathcal{A}}$ is the transformed tensor tubal rank-$k$ approximation of $\mathcal{A}$ defined by the L-Gaussian-Sketch algorithm (i.e., Algorithm \ref{alg:l-trp-sketch} with the Gaussian tensor sketching operator selected). Then
	\begin{equation}\label{eq4.24}
		\begin{split}
			\mathbb{E}{\| \mathcal{A}-\hat{\mathcal{A}}\|} _{F}^{2}\le (1+f(k,s)){(1+\frac{2\varrho}{k-\varrho-1})}\tau_{\varrho+1}^{2}(\mathcal{A}\star{_{L}}\mathcal{A}^{H}),
		\end{split}
	\end{equation}
	where $\varrho$ is a natural number less than $k-1$, $f(\varrho,k):={\varrho}/{(k-\varrho-1)}$, and the tail energy $\tau^{2}_{\varrho+1}$
	is defined by Definition \ref{df2.7}.
\end{theorem}

\textbf{proof}\quad We have
	\begin{align*}
		& \ \mathbb{E}{\| \mathcal{A}-\hat{\mathcal{A}}\|} _{F}^{2}\\
		= & \ \mathbb{E}_{ \Upsilon}\mathbb{E}_{\Omega}{\| \mathcal{A}-\mathcal{Q}\star{_{L}} \mathcal{Q}^{H}\star{_{L}} \mathcal{A}*\mathcal{P}\star{_{L}} \mathcal{P}^{H} \|} _{F}^{2}+\mathbb{E}{\| \mathcal{C}-\mathcal{Q}^{H}\star{_{L}} \mathcal{A}\star{_{L}} \mathcal{P} \|} _{F}^{2}\\
		= & \ (1+f(k,s))\mathbb{E}_{\Upsilon}\mathbb{E}_{\Omega}{\| \mathcal{A}-\mathcal{Q}\star{_{L}} \mathcal{Q}^{H}\star{_{L}} \mathcal{A}\star{_{L}} \mathcal{P}\star{_{L}} \mathcal{P}^{H} \|} _{F}^{2}\\
		& \ + \frac{k(2k+1-s)}{(s-k-1)^{2}}\mathbb{E}{\| \mathcal{Q}_{\perp}^{H}\star{_{L}} \mathcal{A}\star{_{L}} \mathcal{P}_{\perp} \|} _{F}^{2}\\
		\le & \ (1+f(k,s)){(1+\frac{2\varrho}{k-\varrho-1})}\tau_{\varrho+1}^{2}(\mathcal{A}\star{_{L}}\mathcal{A}^{H}) \\
		& \ +\frac{k(2k+1-s)}{(s-k-1)^{2}}\mathbb{E}{\| \mathcal{Q}_{\perp}^{H}\star{_{L}} \mathcal{A}\star{_{L}} \mathcal{P}_{\perp} \|} _{F}^{2}\\
		\le & \ (1+f(k,s)){(1+\frac{2\varrho}{k-\varrho-1})}\tau_{\varrho+1}^{2}(\mathcal{A}\star{_{L}}\mathcal{A}^{H}).
	\end{align*}
	In particular, the first equation is known by Proposition \ref{pro4.3} in Appendix B; the second is known by Lemma \ref{lem4.6} in Appendix B; the first inequality is known by Theorem \ref{thm4.4} in Appendix B; and the last inequality is because we require $s\ge 2k+1$, and thus the missing item $\frac{k(2k+1-s)}{(s-k-1)^{2}}\mathbb{E}{\| \mathcal{Q}_{\perp}^{H}\star{_{L}} \mathcal{A}\star{_{L}} \mathcal{P}_{\perp} \|} _{F}^{2}$ is negative. This completes the proof. 

\begin{theorem}\label{thmPI}
	Assume that the sketch parameters satisfy $s\ge 2k+1$ and $\hat{\mathcal{A}}$ is the transformed tensor tubal rank-$k$ approximation of $\mathcal{A}$ defined by Algorithm \ref{alg:DCT-Sketch-PI}. Then
	\begin{equation}\label{eq4.PI}
		\begin{split}
			\mathbb{E}{\| \mathcal{A}-\hat{\mathcal{A}}\|} _{F}^{2}\le (1+f(k,s)){(1+\frac{2\varrho}{k-\varrho-1})}\tau_{\varrho+1}^{2}{((\mathcal{A}\star{_{L}}\mathcal{A}^{H})^{(2q+1)})},
		\end{split}
	\end{equation}
	where $q$ is a nonnegative integer, $\varrho$ is a natural number less than $k-1$, $f(\varrho,k):={\varrho}/{(k-\varrho-1)}$,  the tail energy $\tau^{2}_{\varrho+1}$
	is defined by Definition \ref{df2.7}, and $L$ here is $C$ in Algorithm \ref{alg:DCT-Sketch-PI}.
\end{theorem}
\textbf{proof}\quad 
	We exploit the third order tensor $(\mathcal{A} \star{_{L}} \mathcal{A}^H)^q \star{_{L}} \mathcal{A}$ instead of the original tensor $\mathcal{A}$. Following the proof in Theorem \ref{thm4.7}, we have
	
	\begin{align*}
		& \ \mathbb{E}{\| \mathcal{A}-\hat{\mathcal{A}}\|} _{F}^{2}\\
		= & \ \mathbb{E}_{ \Upsilon}\mathbb{E}_{\Omega}{\| \mathcal{A}-\mathcal{Q}\star{_{L}} \mathcal{Q}^{H}\star{_{L}} (\mathcal{A} \star{_{L}} \mathcal{A}^H)^q \star{_{L}} \mathcal{A} \star{_{L}}\mathcal{P}\star{_{L}} \mathcal{P}^{H} \|} _{F}^{2}\\
		& \ +\mathbb{E}{\| \mathcal{C}-\mathcal{Q}^{H}\star{_{L}} (\mathcal{A} \star{_{L}} \mathcal{A}^H)^q \star{_{L}} \mathcal{A} \star{_{L}} \mathcal{P} \|} _{F}^{2}\\
		= & \ (1+f(k,s))\mathbb{E}_{\Upsilon}\mathbb{E}_{\Omega}{\| \mathcal{A}-\mathcal{Q}\star{_{L}} \mathcal{Q}^{H}\star{_{L}} (\mathcal{A} \star{_{L}} \mathcal{A}^H)^q \star{_{L}} \mathcal{A} \star{_{L}} \mathcal{P}\star{_{L}} \mathcal{P}^{H} \|} _{F}^{2}\\
		& \ + \frac{k(2k+1-s)}{(s-k-1)^{2}}\mathbb{E}{\| \mathcal{Q}_{\perp}^{H}\star{_{L}} (\mathcal{A} \star{_{L}} \mathcal{A}^H)^q \star{_{L}} \mathcal{A} \star{_{L}} \mathcal{P}_{\perp} \|} _{F}^{2}\\
		\le & \ (1+f(k,s)){(1+\frac{2\varrho}{k-\varrho-1})}\tau_{\varrho+1}^{2}\big((\mathcal{A} \star{_{L}} \mathcal{A}^H)^q \star{_{L}} \mathcal{A} \star{_{L}}((\mathcal{A} \star{_{L}} \mathcal{A}^H)^q \star{_{L}} \mathcal{A})^{H}\big) \\
		& \ +\frac{k(2k+1-s)}{(s-k-1)^{2}}\mathbb{E}{\| \mathcal{Q}_{\perp}^{H}\star{_{L}} (\mathcal{A} \star{_{L}} \mathcal{A}^H)^q \star{_{L}} \mathcal{A} \star{_{L}} \mathcal{P}_{\perp} \|} _{F}^{2}\\
		\le & \ (1+f(k,s)){(1+\frac{2\varrho}{k-\varrho-1})}\tau_{\varrho+1}^{2}((\mathcal{A}\star{_{L}}\mathcal{A}^{H})^{2q+1}).
	\end{align*}
	This completes the proof.

\section{\textbf{Numerical Experiments}} \label{sec:experiments}

This section showcases numerical experiments validating the efficiency and effectiveness of the proposed two-sided sketching algorithms, in comparison with the state-of-the-art algorithms including the T-Sketch algorithm \cite[Algorithm 2]{Qi21}, T-Sketch-PI algorithm ($q=1$) \cite{Qi21}, truncated-t-SVD algorithm \cite{Kilmer11}, and rt-SVD algorithm (i.e., a one-sided randomized algorithm based on t-SVD in \cite[Algorithm 6]{Zhang18}). The sketch size parameters are set to $s = 2k+1$. The following relative error $\epsilon_{\rm err}$ and the peak signal-to-noise ratio (PSNR) $\rho_{\rm psnr}$ are used as metrics of the low-rank approximation $\hat{\mathcal{A}}$ to the input tensor data $\mathcal{A}$, i.e.,
\begin{equation}	
	\epsilon_{\rm err} :={{||\mathcal{A}-\hat{\mathcal{A}}||}_{F}^{2}}/{||\mathcal{A}||_{F}^{2}},
	\quad \rho_{\rm psnr}:=10\log_{10}\frac{n_{1}n_{2}n_{3}||\mathcal{A}||_{\infty}^{2}}{||\mathcal{A}-\hat{\mathcal{A}}||_{F}^{2}}.
\end{equation}


\subsection{Synthetic Experiments}
We firstly conduct numerical tests on some synthetic input tensors $\mathcal{A}\in \mathbb{R}^{10^{3}\times 10^{3}\times 10}$ with decaying spectrum.


{\it Polynomial decay:}	These tensors are f-diagonal tensors. Considering their $j$-th frontal slices with the form 
\begin{equation}
	\mathcal{A}^{(j)}=\mbox{diag}\big(\underbrace{1,\ldots,1}_{\min(r,j)},\\2^{-p},3^{-p},4^{-p},\ldots,(n-\min(r,j)+1)^{-p}\big) \in \mathbb{R}^{n\times n},
\end{equation}
we study two examples, i.e., PolyDecaySlow ($p = 0.5$) 
and PolyDecayFast ($p = 2$).

\begin{figure}
	\centering
	\includegraphics[trim={{.30\linewidth} {.08\linewidth} {.20\linewidth} {.10\linewidth}}, clip, width=0.87\linewidth, height = 0.47\linewidth]{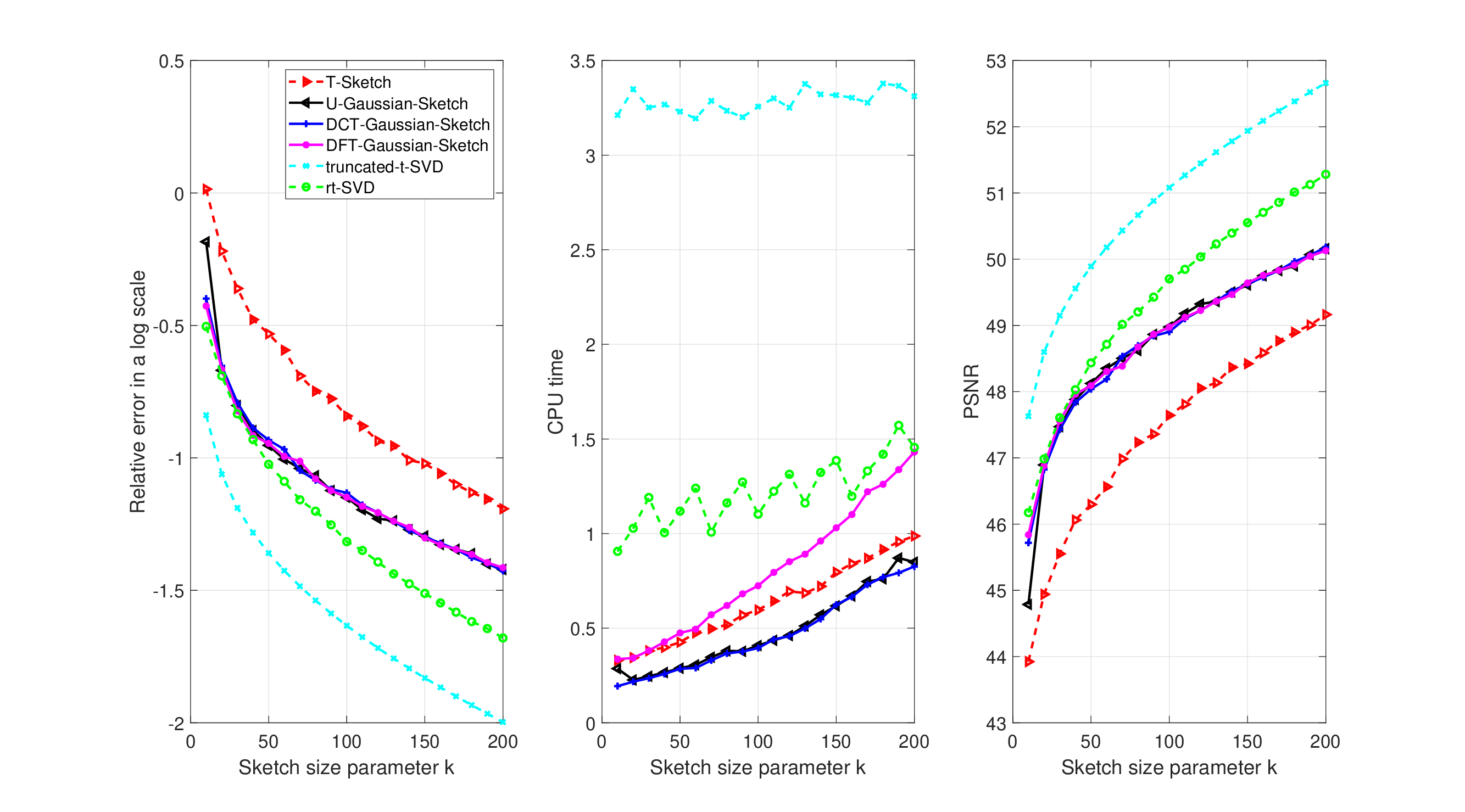} \vspace{0.05in} \\
	\includegraphics[trim={{.30\linewidth} {.08\linewidth} {.20\linewidth} {.10\linewidth}}, clip, width=0.87\linewidth, height = 0.47\linewidth]{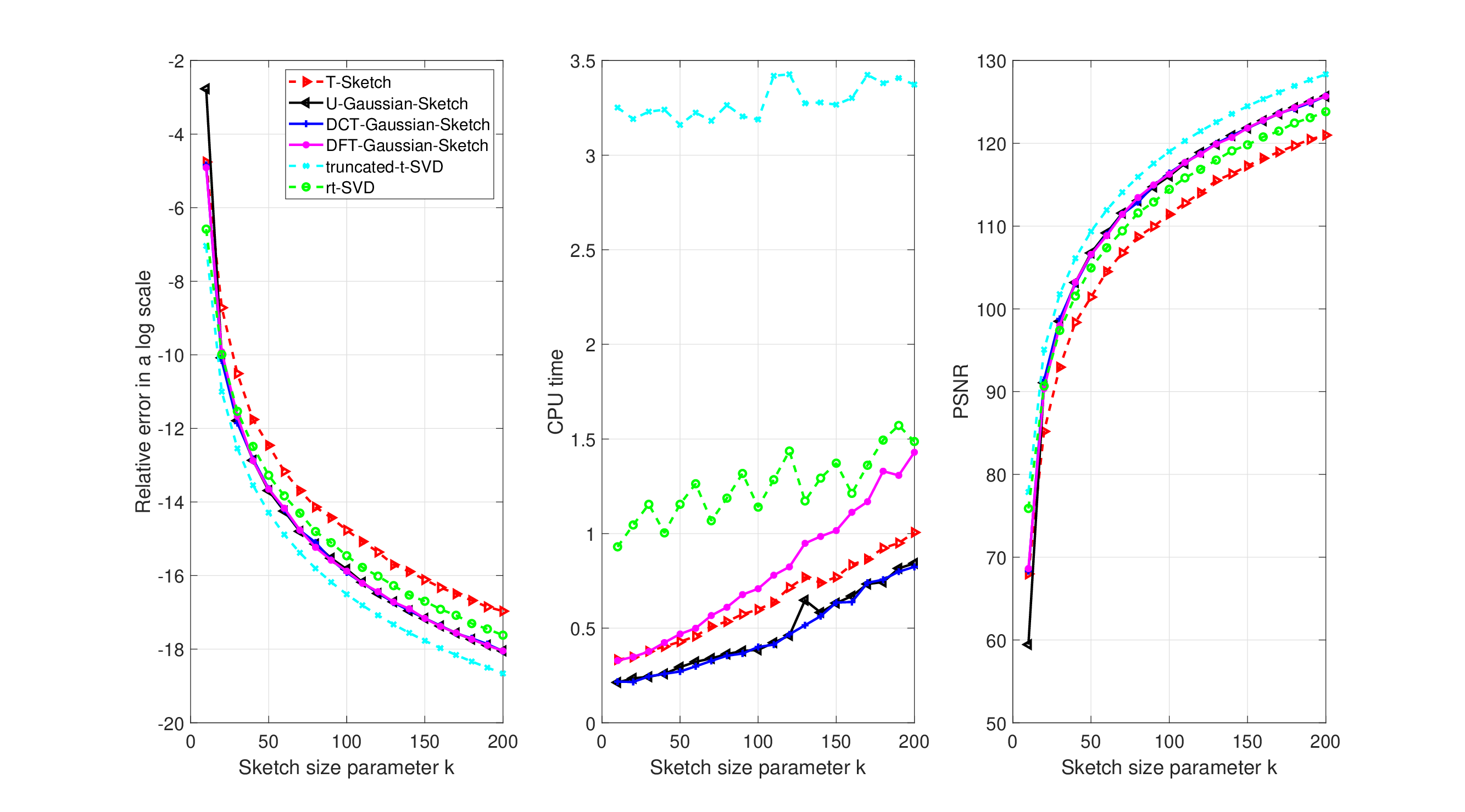}
	\caption{Performance of different methods and our method via different transformed domains (i.e., U, DCT, and DFT transforms) in terms of the relative error ({\it left} column), CPU time ({\it middle} column), and PSNR ({\it right} column). In particular, the first row is for tensor data (PolyDecaySlow) with slow decaying spectrum and the second row is for tensor data (PolyDecayFast) with fast decaying spectrum.}
\label{fig:syn-about-TD}
\end{figure}

\begin{figure}
\centering
\includegraphics[trim={{.30\linewidth} {.08\linewidth} {.20\linewidth} {.10\linewidth}}, clip, width=0.87\linewidth, height = 0.47\linewidth]{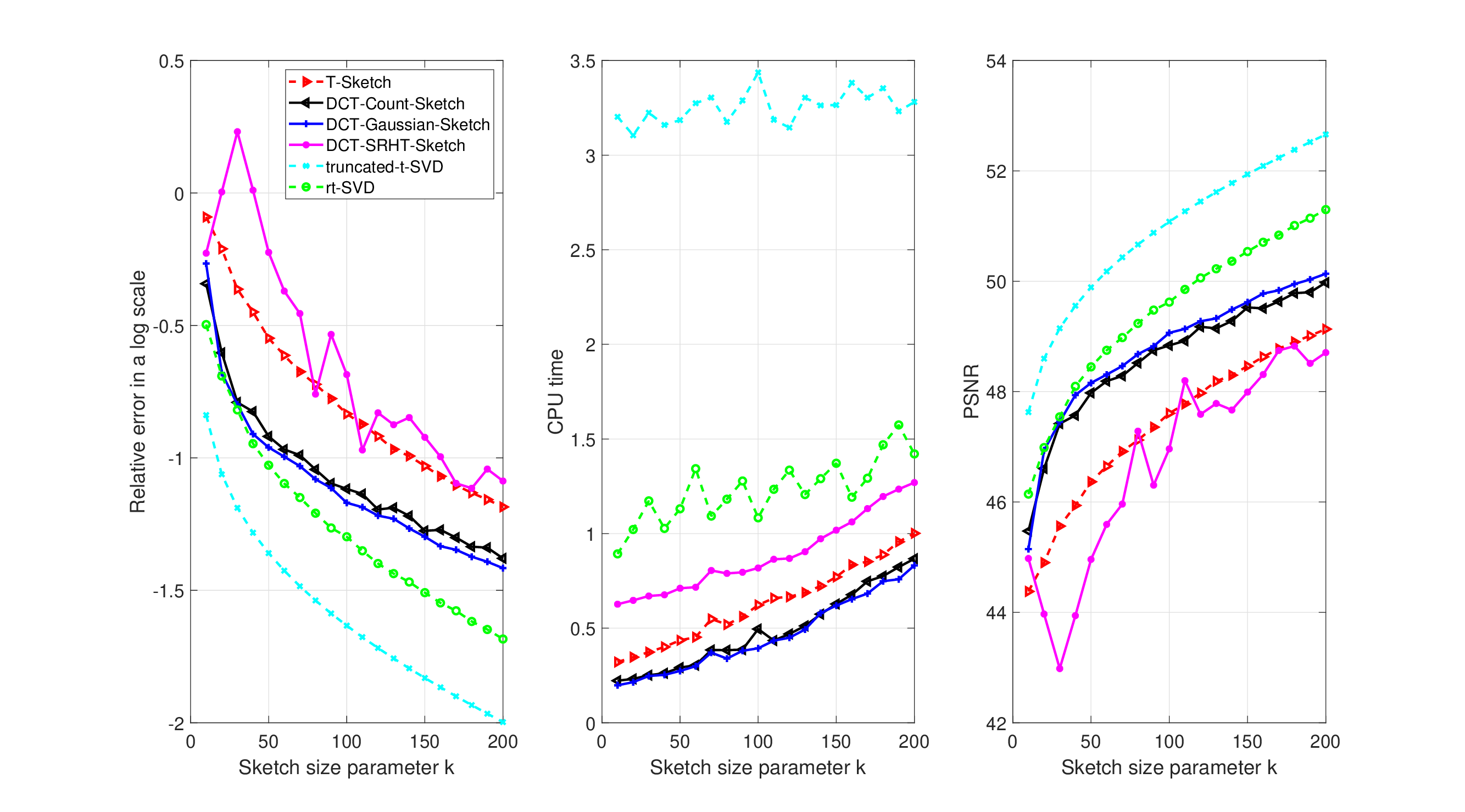}  \vspace{0.05in} \\
\includegraphics[trim={{.30\linewidth} {.08\linewidth} {.20\linewidth} {.10\linewidth}}, clip, width=0.87\linewidth, height = 0.47\linewidth]{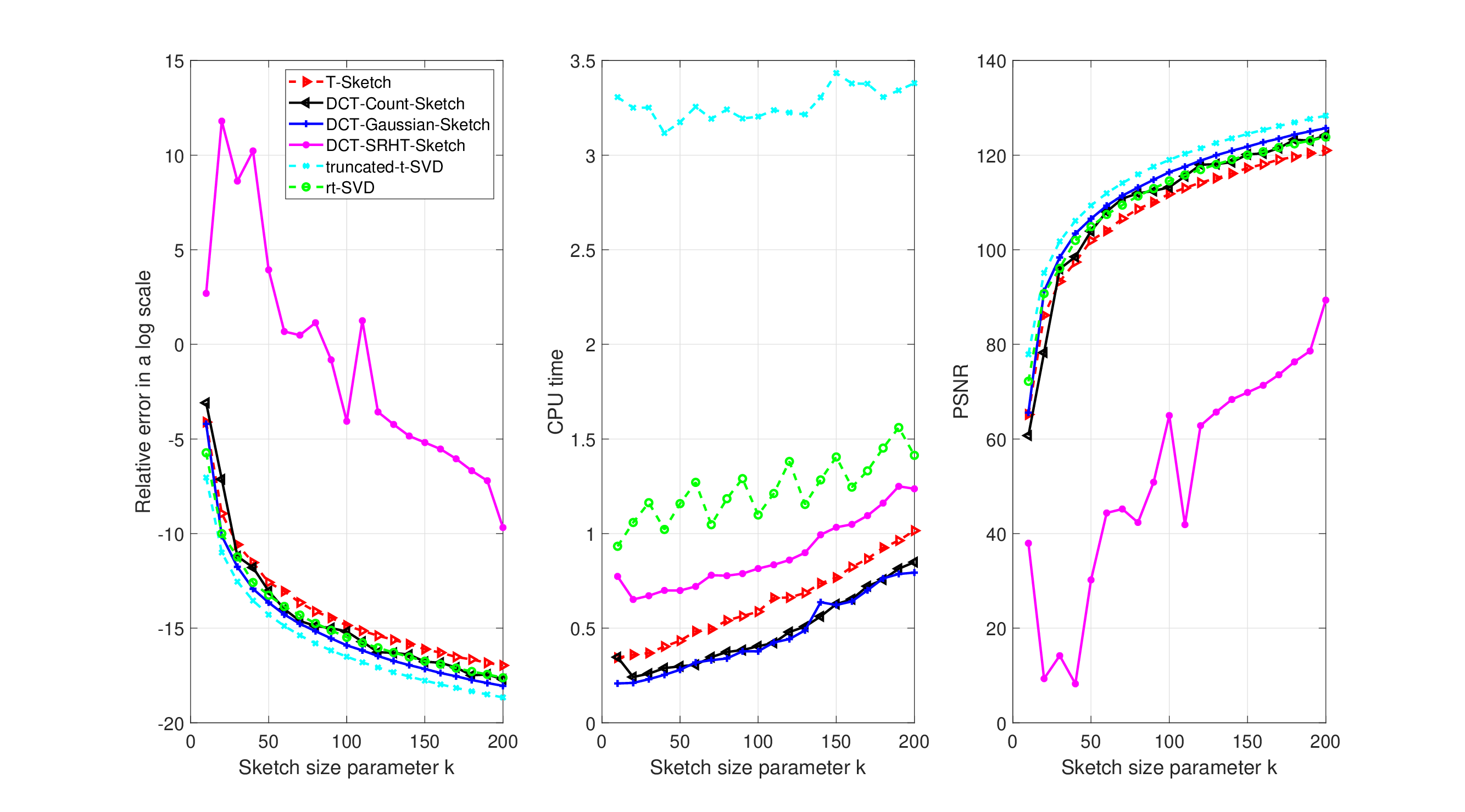}
\caption{Performance of different methods and our method based on the DCT transformed domain and via different tensor sketching operators (i.e., Gaussian, SRHT, and count sketching operators) in terms of the relative error ({\it left} column), CPU time ({\it middle} column), and PSNR ({\it right} column). In particular, the first row is for tensor data (PolyDecaySlow) with slow decaying spectrum and the second row is for tensor data (PolyDecayFast) with fast decaying spectrum.}
\label{fig:syn-about-TSO}
\end{figure}

\begin{figure}
\centering
\includegraphics[trim={{.30\linewidth} {.08\linewidth} {.20\linewidth} {.10\linewidth}}, clip, width=0.87\linewidth, height = 0.47\linewidth]{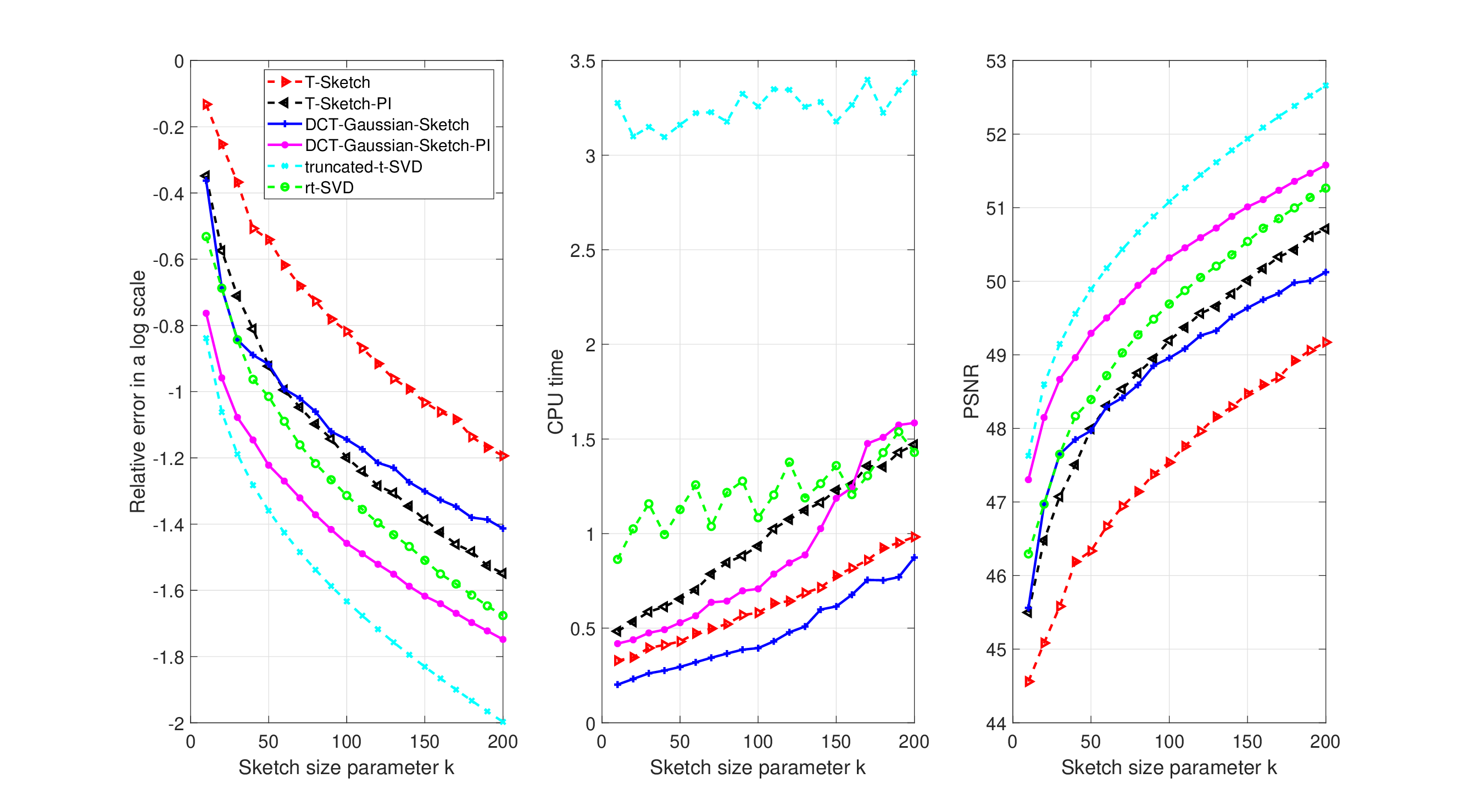} 
\vspace{0.05in} \\
\includegraphics[trim={{.30\linewidth} {.08\linewidth} {.20\linewidth} {.10\linewidth}}, clip, width=0.87\linewidth, height = 0.47\linewidth]{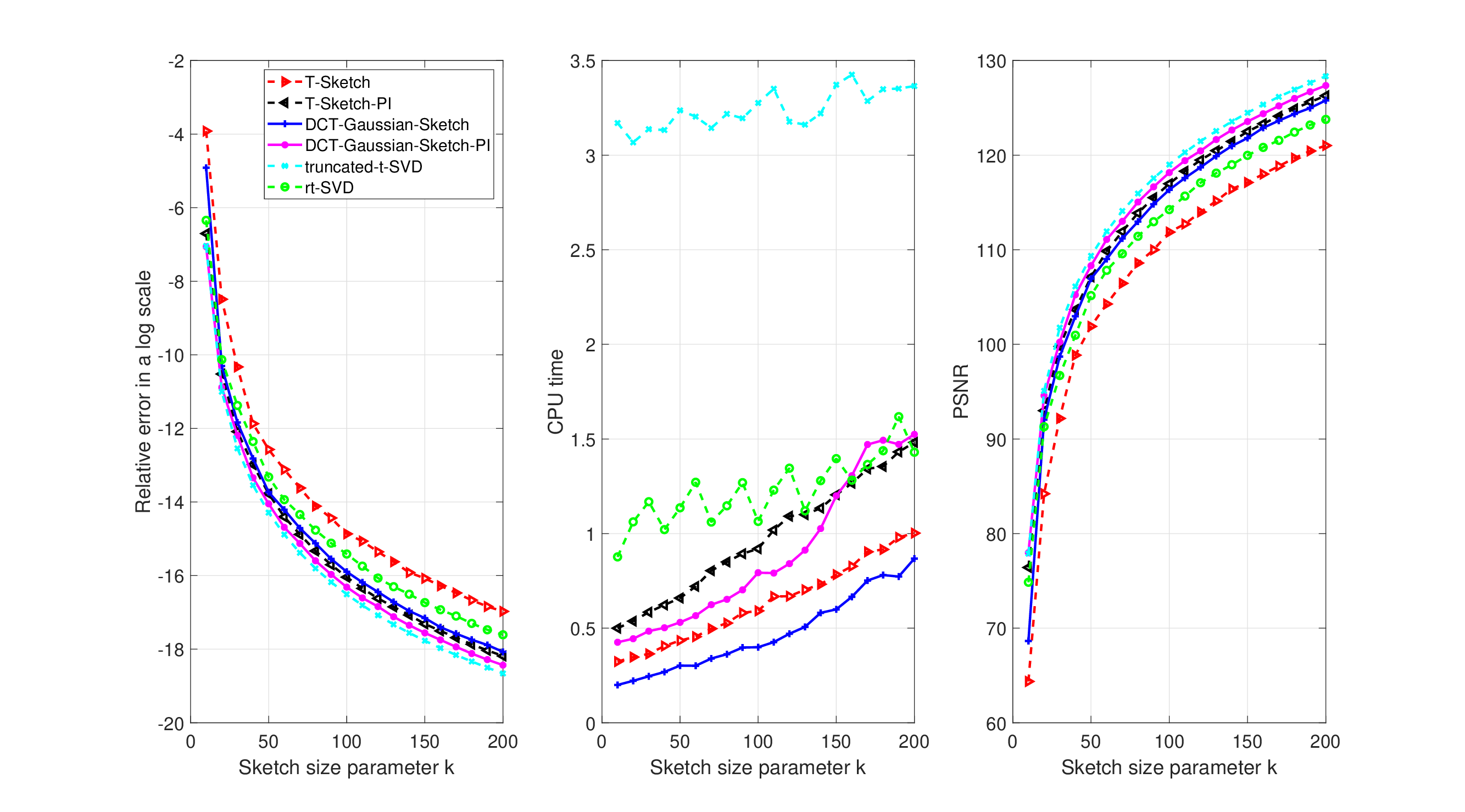}
\caption{Performance of different methods and our method based on the DCT transformed domain and the Gaussian tensor sketching operator with and without the power iteration technique in terms of the relative error ({\it left} column), CPU time ({\it middle} column), and PSNR ({\it right} column). In particular, the first row is for tensor data (PolyDecaySlow) with slow decaying spectrum and the second row is for tensor data (PolyDecayFast) with fast decaying spectrum.
}
\label{fig:syn-about-PI}
\end{figure}

Figures \ref{fig:syn-about-TD}--\ref{fig:syn-about-PI} give the results of the algorithms compared in terms of the relative error, CPU time, and PSNR, as the size of the sketch parameter $k$ varies; see more details below.

\subsubsection{Experiments Regarding the Transformed Domains}
We first examine the performance of different transformed domains in our two-sided sketching method.
Figure \ref{fig:syn-about-TD} shows the results of different methods and our method via different transformed domains (specifically, U-Gaussian-Sketch, DCT-Gaussian-Sketch, and DFT-Gaussian-Sketch algorithms) in terms of the relative error, CPU time, and PSNR.
Figure \ref{fig:syn-about-TD} illustrates that among our two-sided Gaussian sketching algorithms with different transformed domains, they all perform similarly in terms of the relative error and PSNR. Regarding the CPU time, our DCT-Gaussian-Sketch and U-Gaussian-Sketch algorithms outperform all other methods including the rt-SVD, truncated-t-SVD, and T-Sketch algorithms. This means, regarding the two-sided Gaussian sketching algorithms employing various transformation operators, the DCT and U transforms emerge as the most efficient. 

As shown in the second row of Figure \ref{fig:syn-about-TD}, our two-sided sketching method with different transformed domains all  surpasses the accuracy of the rt-SVD and T-Sketch algorithms for input tensors exhibiting a fast decay spectrum. Consequently, with reduced storage requirements and manipulation, our two-sided sketching method with different transformed domains provides superior accuracy in low-rank approximations. It is worth highlighting that, for input tensors exhibiting a fast decay spectrum, our method is the second only to the truncated-t-SVD algorithm in terms of accuracy, but our method offers a significant speed advantage, i.e, the truncated-t-SVD is far slower than our method (see the middle plot of the second column of Figure \ref{fig:syn-about-TD}).

\subsubsection{Experiments Regarding the Tensor Sketching Operators}
We now examine the performance of different tensor sketching operators in our two-sided sketching method. As an example, our method here is adopted with the DCT transformed domain.
Figure \ref{fig:syn-about-TSO} shows the results of different methods and our method via different tensor sketching operators (specifically, DCT-Gaussian-Sketch, DCT-SRHT-Sketch, DCT-Count-Sketch algorithms) in terms of the relative error, CPU time, and PSNR.
Figure \ref{fig:syn-about-TSO} illustrates that among our two-sided sketching algorithms based on the DCT transformed domain, the DCT-Gaussian-Sketch algorithm outperforms the other two in terms of accuracy, indicating the better effectiveness of the Gaussian tensor sketching operator.
Regarding the CPU time, our DCT-Gaussian-Sketch and DCT-Count-Sketch algorithms outperform all other methods including the rt-SVD, truncated-t-SVD, and T-Sketch algorithms. This means, regarding the two-sided sketching algorithms (based on the DCT transformed domain) employing various tensor sketching operators, the Gaussian and count sketching operators emerge as the most efficient. Considering both the accuracy and speed, the DCT-Gaussian-Sketch algorithm is the best among our method with different tensor sketching operators.

As shown in the second row of Figure \ref{fig:syn-about-TSO}, consistent results are obtained for our DCT-Gaussian-Sketch algorithm as that in Figure \ref{fig:syn-about-TD}. It is worth highlighting that, for input tensors exhibiting a fast decay spectrum, our DCT-Gaussian-Sketch algorithm is the second only to the truncated-t-SVD algorithm in terms of accuracy, but it offers a significant speed advantage, i.e., the truncated-t-SVD is far slower than our DCT-Gaussian-Sketch algorithm (see the middle plot of the second column of Figure~\ref{fig:syn-about-TSO}).

\begin{figure}[htb]
\centering
\begin{tabular}{cc}
	\includegraphics[trim={{5.0in} {6.5in} {19.0in} {5.5in}}, clip, width=2.0in, height = 1.3in]{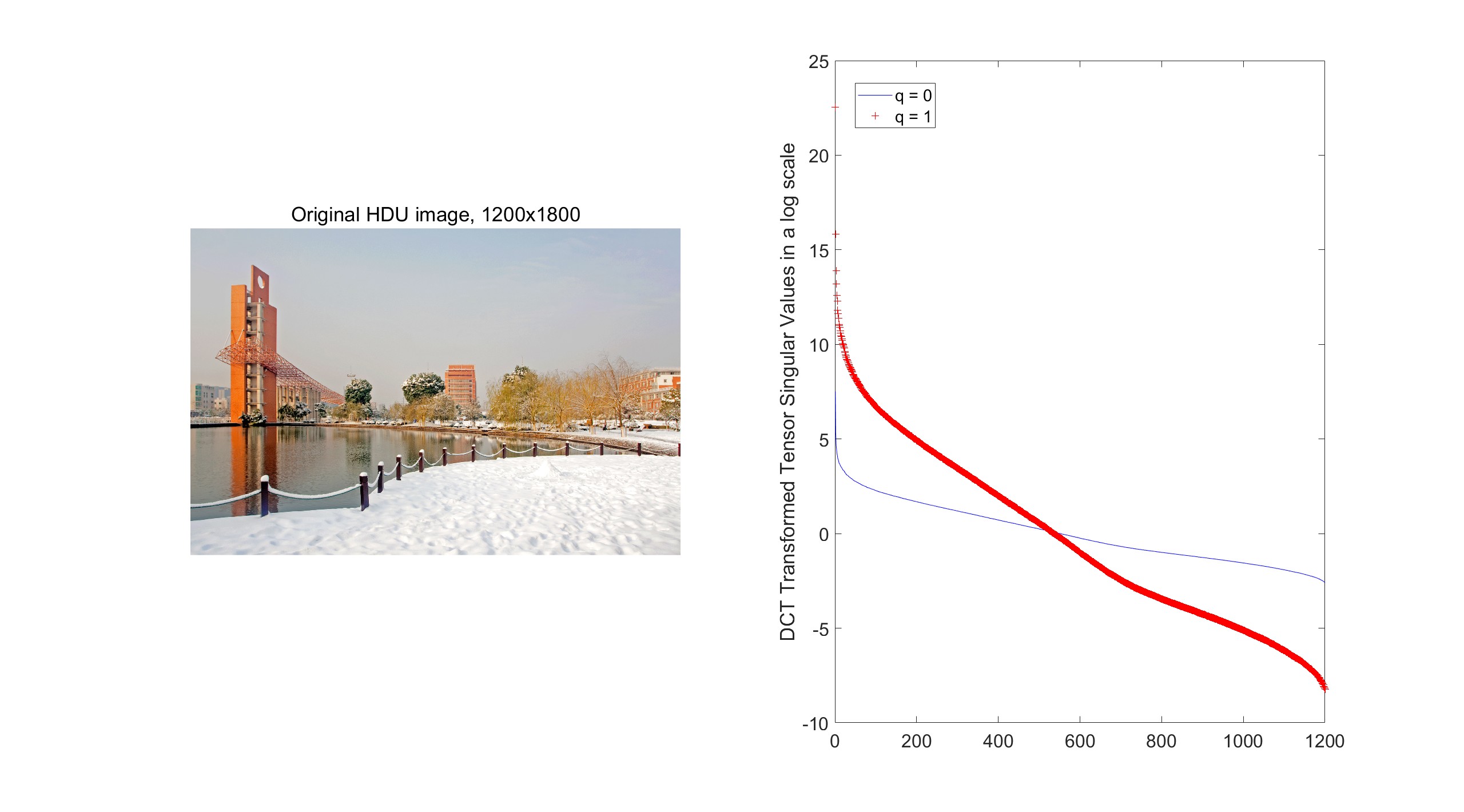} &
	\includegraphics[trim={{5.0in} {6.5in} {19.0in} {5.5in}}, clip, width=2.0in, height = 1.3in]{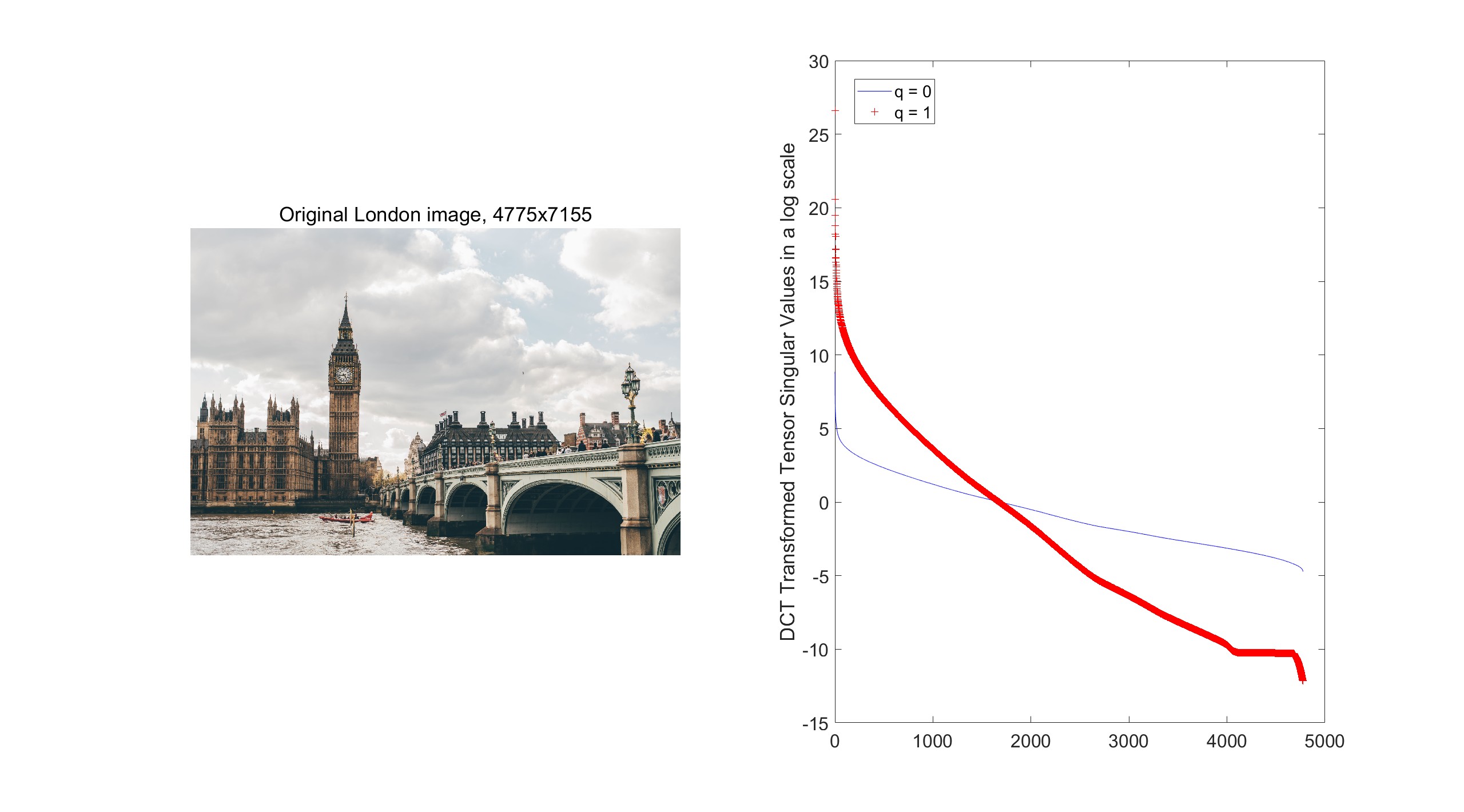} \\
	{\footnotesize (a) HDU image (size: $1,200\times 1,800$)} & {\footnotesize (b)  London image (size: $4,775\times 7,155$)} \\
	\includegraphics[trim={{18.5in} {1.5in} {2.9in} {1.0in}}, clip, width=2.2in, height = 3.0in]{HDUSingular.jpg} &
	\includegraphics[trim={{18.5in} {1.5in} {2.9in} {1.0in}}, clip, width=2.2in, height = 3.0in]{LONDONSingular.jpg} \\
	{\footnotesize (c) Transformed tensor singular values of (a)} & {\footnotesize (d) Transformed tensor singular values of (b)} 
\end{tabular}
\caption{Test real-world color images and the DCT transformed tensor singular values in terms of their decaying spectrum. Row one: the original color images. Row two: the transformed tensor singular values of the power iteration ($q=0$) and ($q=1$) for the given color images.}
\label{fig:image-decay-spec}
\end{figure}

\subsubsection{Experiments Regarding the Power Iteration Technique}
Figure~\ref{fig:syn-about-PI} shows the results of different methods and our method based on the DCT transformed domain and the Gaussian tensor sketching operator with and without the power iteration technique (i.e., DCT-Gaussian-Sketch and DCT-Gaussian-Sketch-PI ($q=1$)) in terms of the relative error, CPU time, and PSNR. For the input tensors exhibiting both slow and fast decay spectra, Figure \ref{fig:syn-about-PI} shows that DCT-Gaussian-Sketch-PI can indeed outperform DCT-Gaussian-Sketch in terms of accuracy, with a litter sacrifice of speed, indicating the effectiveness of the power iteration technique in approximation quality enhancement. Moreover, for the comparison between our DCT-Gaussian-Sketch-PI algorithm, the T-Sketch algorithm with the power iteration technique (i.e., T-Sketch-PI), and the rt-SVD algorithm, our method is superior in all metrics. 

With reduced storage requirements and manipulation, as shown in the second row of Figure \ref{fig:syn-about-PI} for input tensors exhibiting a fast decay spectrum, 
our DCT-Gaussian-Sketch-PI algorithm is the second only to the truncated-t-SVD algorithm in terms of accuracy, but our algorithm again offers a significant speed advantage, i.e, the truncated-t-SVD is far slower than our DCT-Gaussian-Sketch-PI algorithm (see the middle plot of the second column of Figure \ref{fig:syn-about-PI}).



\begin{figure}
\centering
\includegraphics[trim={{.30\linewidth} {.07\linewidth} {.20\linewidth} {.06\linewidth}}, clip, width=0.87\linewidth, height = 0.47\linewidth]{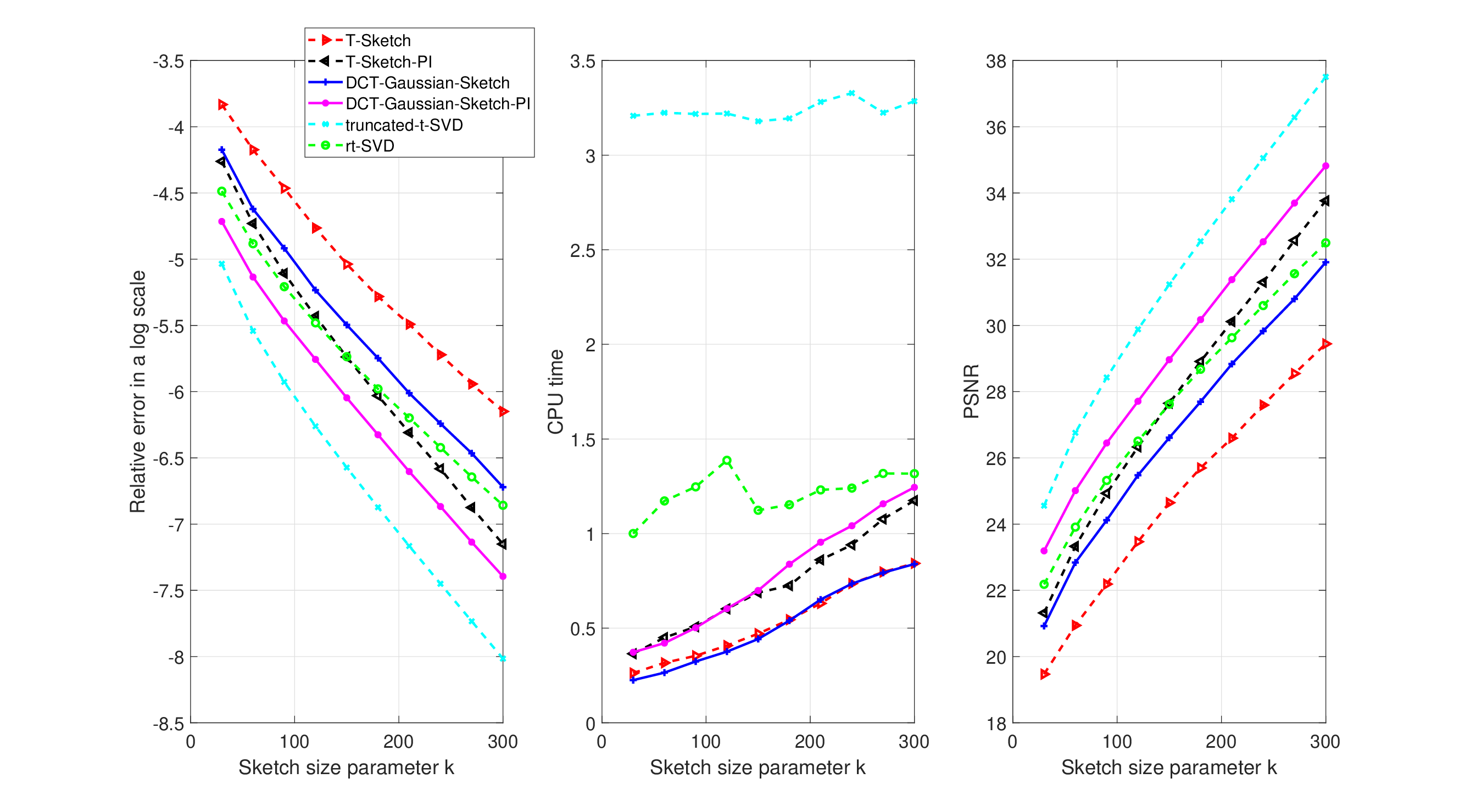} 
\vspace{0.05in} \\
\includegraphics[trim={{.30\linewidth} {.07\linewidth} {.20\linewidth} {.08\linewidth}}, clip, width=0.87\linewidth, height = 0.47\linewidth]{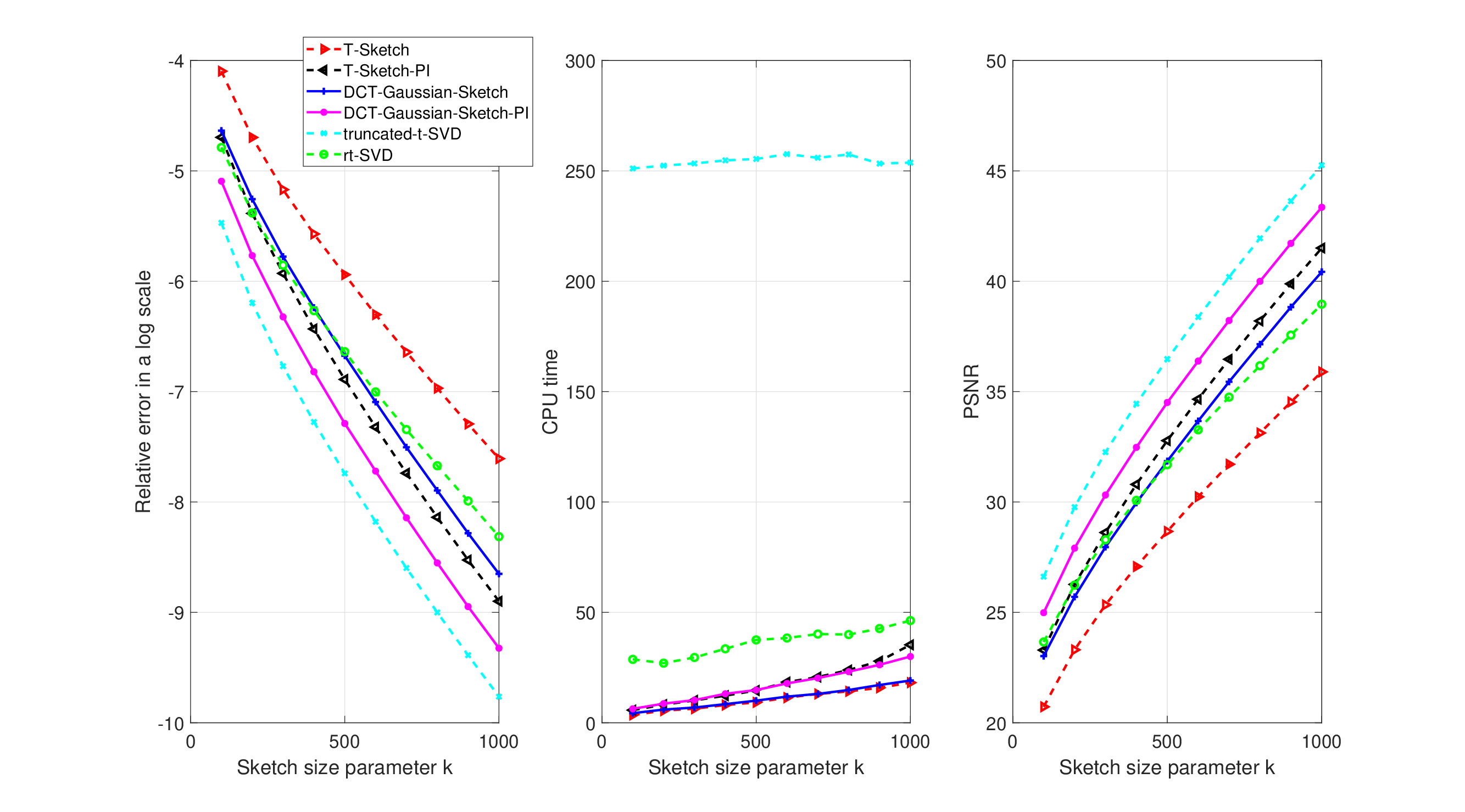}
\caption{
	Low-rank approximation performance of different methods on the HDU image ({\it first} row) and the London image ({\it second} row) in terms of the relative error ({\it left} column), CPU time ({\it middle} column), and PSNR ({\it right} column).}
\label{fig:LRA-image-quan}
\end{figure}

\subsection{Real-world Data}
We now conduct experiments on real-world data including color images and grayscale videos. For our method, we select the DCT-Gaussian-Sketch and DCT-Gaussian-Sketch-PI ($q=1$) algorithms, given their great performance demonstrated in the previous section.

\subsubsection{Color Images}
Two large size color images, i.e., HDU picture\footnote{https://www.hdu.edu.cn/landscape} with size of $1,200\times1,800\times3$ and London picture with size of $4,775\times7,155\times3$, are employed, see the first row of Figure \ref{fig:image-decay-spec}.  

The second row of Figure \ref{fig:image-decay-spec} gives the ordered transformed tensor singular values (indicating the decaying spectrum) using the DCT transformed tensor singular values  decomposition on the two original color images. It clearly shows that the one with power iteration achieves a faster decaying spectrum compared to that without power iteration, demonstrating the great effectiveness of the power iteration technique.


%

\begin{figure}[htbp!]
\centering			       \begin{tabular}{ccc}
	\includegraphics[trim={{4.6in} {12.5in} {23.2in} {2.3in}}, clip, width=1.5in, height = 1.0in]{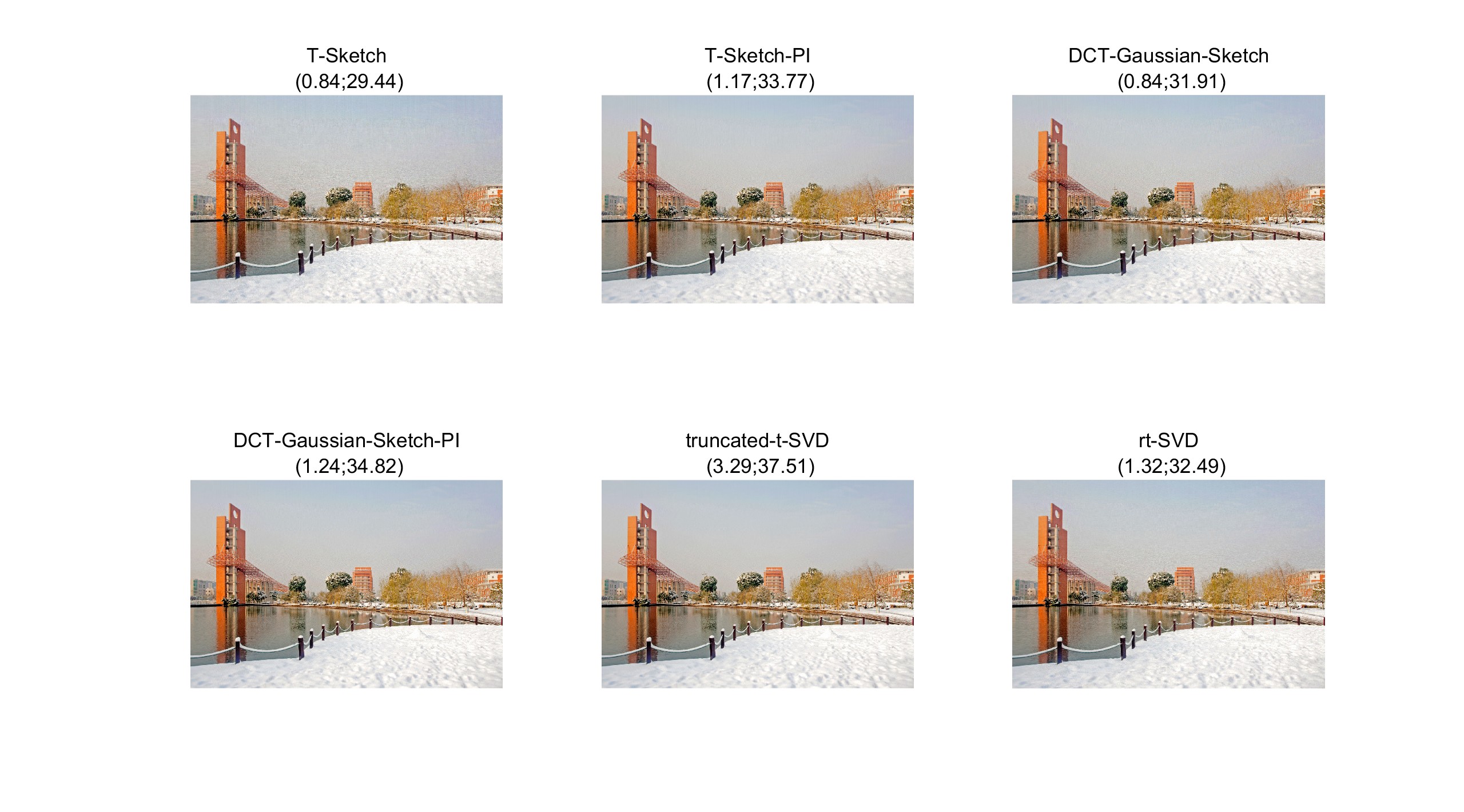} &
	\includegraphics[trim={{14.6in} {12.5in} {13.4in} {2.3in}}, clip, width=1.5in, height = 1.0in]{GaussianPIHDUK=300.jpg} &
	\includegraphics[trim={{24.5in} {12.5in} {3.3in} {2.3in}}, clip, width=1.5in, height = 1.0in]{GaussianPIHDUK=300.jpg} 
	\vspace{-0.03in} \\
	{\footnotesize (a1) T-Sketch} & {\footnotesize (b1) T-Sketch-PI} & {\footnotesize (c1) DCT-Gaussian-Sketch }  \vspace{-0.03in} \\
	{\footnotesize CPU:0.84; PSNR: 29.44}	& {\footnotesize CPU:1.17; PSNR: 33.77} & {\footnotesize CPU: 0.84; PSNR: 31.91} \\
	\includegraphics[trim={{4.6in} {3.0in} {23.2in} {11.7in}}, clip, width=1.5in, height = 1.0in]{GaussianPIHDUK=300.jpg} &
	\includegraphics[trim={{14.6in} {3.0in} {13.4in} {11.7in}}, clip, width=1.5in, height = 1.0in]{GaussianPIHDUK=300.jpg} &
	\includegraphics[trim={{24.5in} {3.0in} {3.3in} {11.7in}}, clip, width=1.5in, height = 1.0in]{GaussianPIHDUK=300.jpg}  
	\vspace{-0.03in}\\
	{\footnotesize (d1) DCT-Gaussian-Sketch-PI} & {\footnotesize (e1) truncated-t-SVD} & {\footnotesize (f1) rt-SVD}  \vspace{-0.03in}  \\
	{\footnotesize CPU: 1.24; PSNR: 34.82} & {\footnotesize CPU: 3.29; PSNR: 37.51} & {\footnotesize CPU: 1.32; PSNR: 32.49} \\ 
	
	\includegraphics[trim={{4.6in} {12.5in} {23.2in} {2.3in}}, clip, width=1.5in, height = 1.0in]{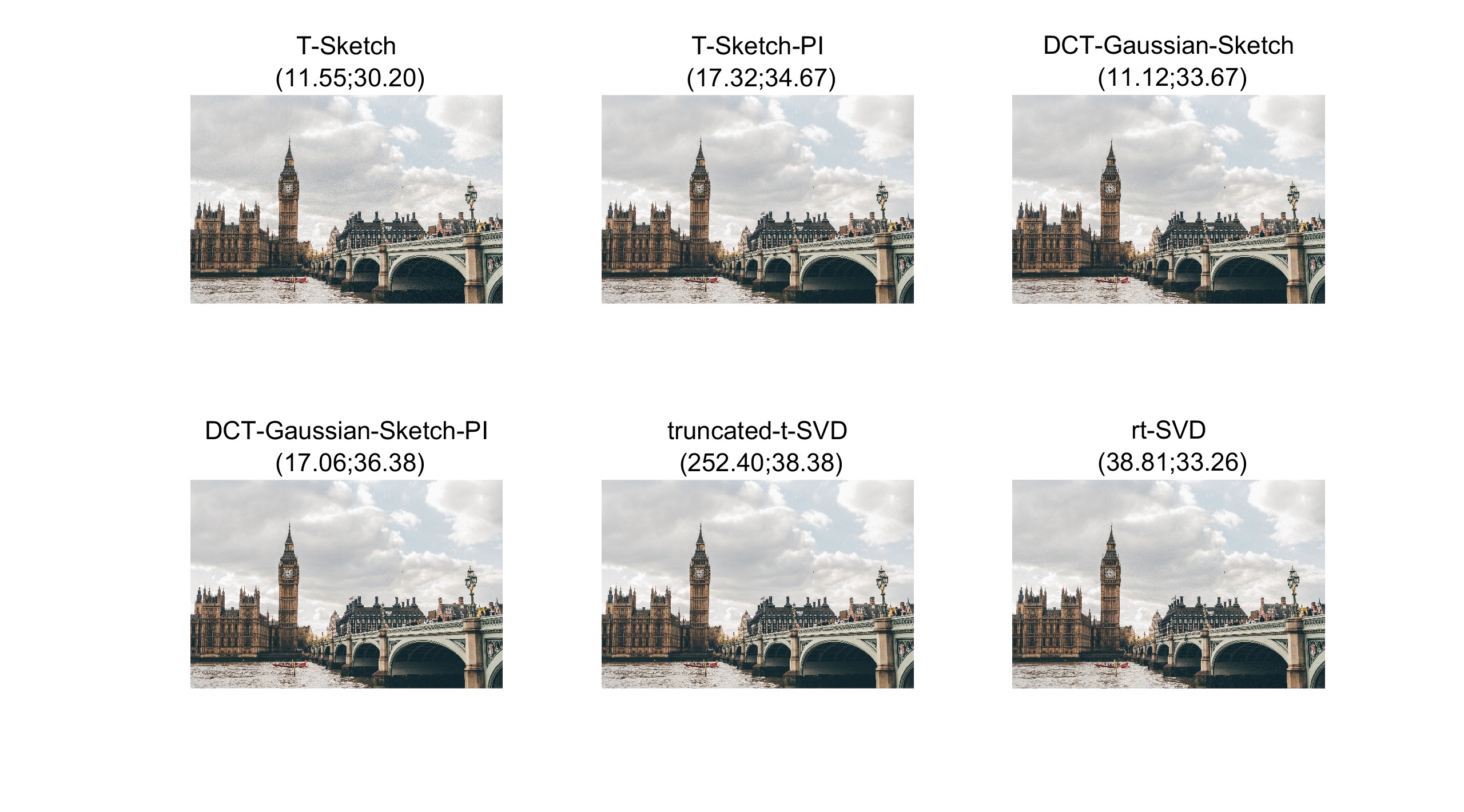} &
	\includegraphics[trim={{14.6in} {12.5in} {13.4in} {2.3in}}, clip, width=1.5in, height = 1.0in]{GaussianPILONDONK=600.jpg} &
	\includegraphics[trim={{24.5in} {12.5in} {3.3in} {2.3in}}, clip, width=1.5in, height = 1.0in]{GaussianPILONDONK=600.jpg} 
	\vspace{-0.03in} \\
	{\footnotesize (a2) T-Sketch} & {\footnotesize (b2) T-Sketch-PI} & {\footnotesize (c2) DCT-Gaussian-Sketch }  \vspace{-0.03in} \\
	{\footnotesize CPU:11.55; PSNR: 30.20}	& {\footnotesize CPU:17.32; PSNR: 34.67} & {\footnotesize CPU: 11.12; PSNR: 33.67} \\
	\includegraphics[trim={{4.6in} {3.0in} {23.2in} {11.7in}}, clip, width=1.5in, height = 1.0in]{GaussianPILONDONK=600.jpg} &
	\includegraphics[trim={{14.6in} {3.0in} {13.4in} {11.7in}}, clip, width=1.5in, height = 1.0in]{GaussianPILONDONK=600.jpg} &
	\includegraphics[trim={{24.5in} {3.0in} {3.3in} {11.7in}}, clip, width=1.5in, height = 1.0in]{GaussianPILONDONK=600.jpg}  
	\vspace{-0.03in} \\
	{\footnotesize (d2) DCT-Gaussian-Sketch-PI} & {\footnotesize (e2) truncated-t-SVD} & {\footnotesize (f2) rt-SVD}  \vspace{-0.03in}  \\
	{\footnotesize CPU: 17.06; PSNR: 36.38} & {\footnotesize CPU: 252.40; PSNR: 38.38} & {\footnotesize CPU: 38.81; PSNR: 33.26} 
\end{tabular}
\caption{Qualitative results of the low-rank approximation performance of different methods on the HDU image (rows 1 and 2) with the sketch size $k$ = 300 and the London image (rows 3 and 4) with the sketch size $k$ = 600 in terms of the CPU time and PSNR.}
\label{fig:LRA-image-quali}	
\end{figure}

Figures \ref{fig:LRA-image-quan} and \ref{fig:LRA-image-quali} respectively give the quantitative and qualitative results of the low-rank approximation performance of different methods on the test color images in terms of different metrics. It is evident that as the sketch parameter size $k$ increases, the low-rank approximation performance of all the methods improves. Consistent results are obtained as that obtained from the previous section for synthetic tensors. Regarding the CPU time, our DCT-Gaussian-Sketch algorithm and the T-Sketch algorithm outperform all others. 
In particular, the second row of Figure \ref{fig:LRA-image-quan} (results for the London image) shows that the DCT-Gaussian-Sketch algorithm surpasses both the rt-SVD and T-Sketch algorithms in terms of  accuracy, with reduced storage and manipulation requirements. The DCT-Gaussian-Sketch-PI algorithm achieves better accuracy  compared to the T-Sketch-PI algorithm with similar speed. The DCT-Gaussian-Sketch-PI algorithm is only second to the truncated-t-SVD algorithm in accuracy, but it is remarkably faster than the truncated-t-SVD algorithm.

\subsubsection{Color Image with Gaussian White Noise}


\begin{figure}[H]
\centering
\includegraphics[trim={{.30\linewidth} {.08\linewidth} {.20\linewidth} {.12\linewidth}}, clip, width=0.87\linewidth, height = 0.47\linewidth]{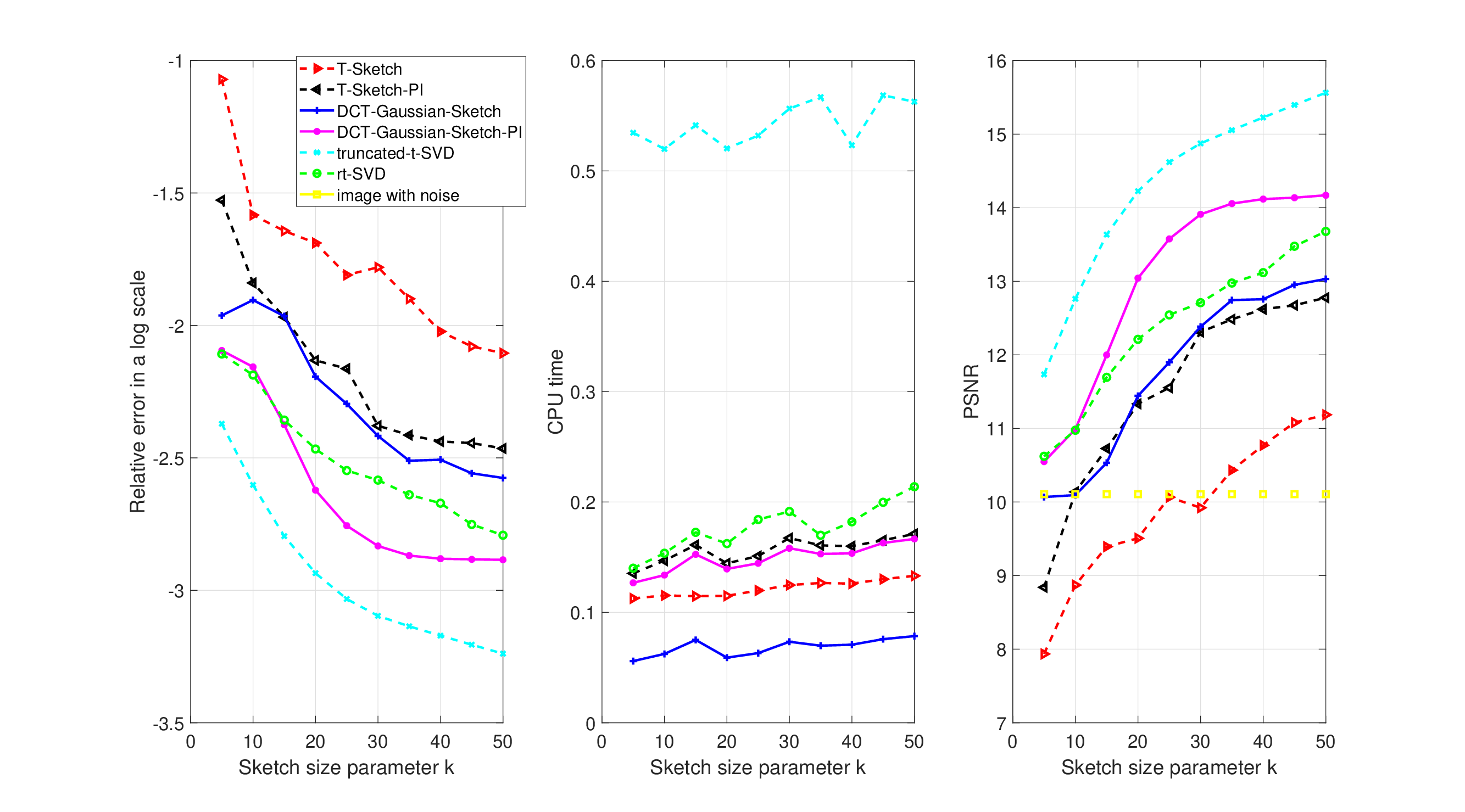} 
\vspace{0.05in} 
\caption{Low-rank approximation performance of different methods on the HDU 2D code image with Gaussian white noise in terms of the relative error ({\it left} column), CPU time ({\it middle} column), and PSNR ({\it right} column). In particular, the yellow line named `image with noise' is utilized to plot the PSNR of the HDU 2D code image with Gaussian white noise.}
\label{fig:WhiteNoise-2D-decay-spec}
\end{figure}

\begin{figure}
\centering			       \begin{tabular}{ccc}
	\includegraphics[trim={{4.6in} {12.5in} {23.2in} {2.3in}}, clip, width=1.5in, height = 1.0in]{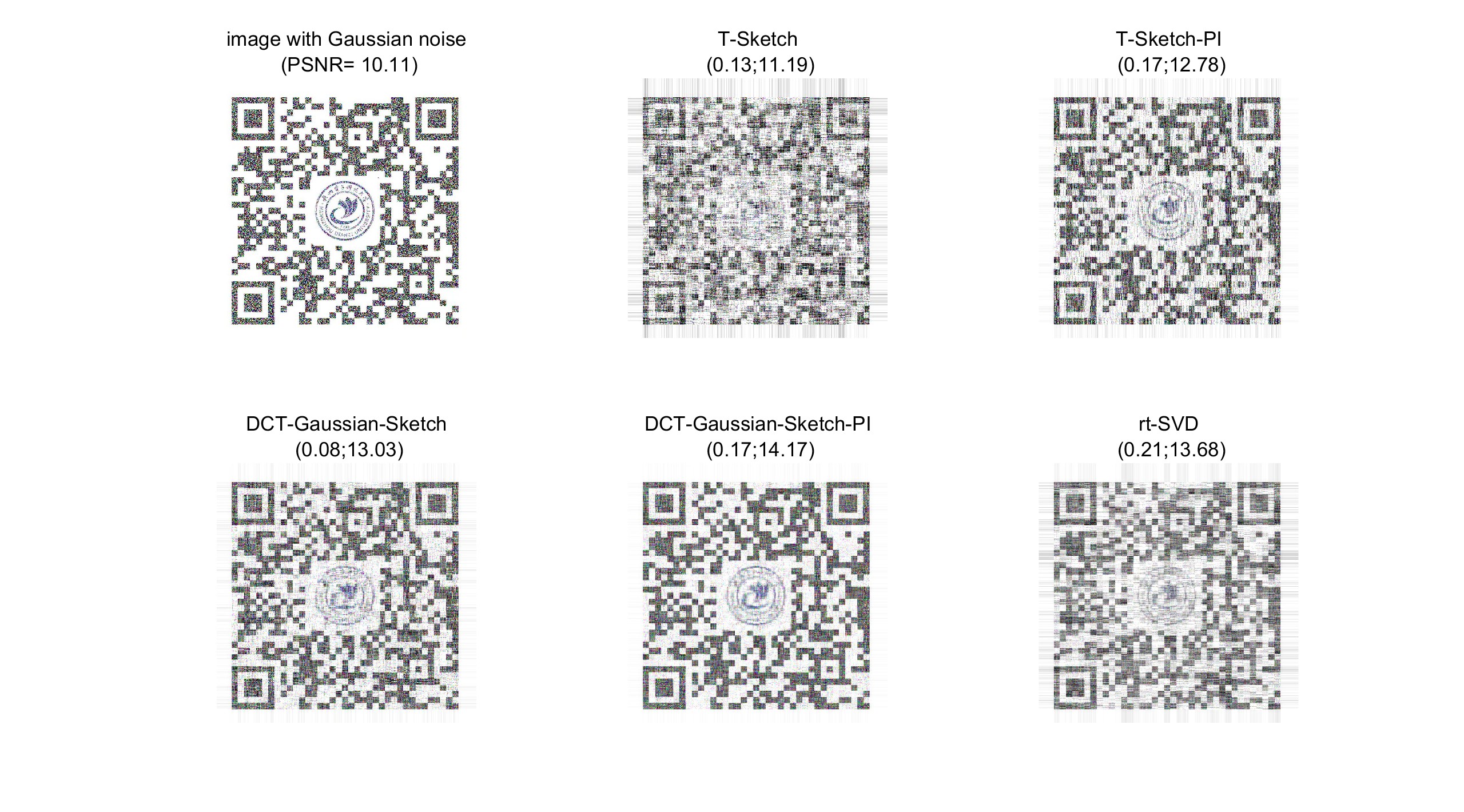} &
	\includegraphics[trim={{14.6in} {12.5in} {13.4in} {2.3in}}, clip, width=1.5in, height = 1.0in]{Noise5ErWeiMaGaussianPIK=50.jpg} &
	\includegraphics[trim={{24.5in} {12.5in} {3.3in} {2.3in}}, clip, width=1.5in, height = 1.0in]{Noise5ErWeiMaGaussianPIK=50.jpg} 
	\vspace{-0.03in} \\
	{\footnotesize (a) Image with noise} & {\footnotesize (b) T-Sketch} & {\footnotesize (c) T-Sketch-PI }  \vspace{-0.03in} \\
	{\footnotesize PSNR: 10.11}	& {\footnotesize CPU:0.13; PSNR: 11.19} & {\footnotesize CPU: 0.17; PSNR: 12.78} \\
	\includegraphics[trim={{4.6in} {3.0in} {23.2in} {11.7in}}, clip, width=1.5in, height = 1.0in]{Noise5ErWeiMaGaussianPIK=50.jpg} &
	\includegraphics[trim={{14.6in} {3.0in} {13.4in} {11.7in}}, clip, width=1.5in, height = 1.0in]{Noise5ErWeiMaGaussianPIK=50.jpg} &
	\includegraphics[trim={{24.5in} {3.0in} {3.3in} {11.7in}}, clip, width=1.5in, height = 1.0in]{Noise5ErWeiMaGaussianPIK=50.jpg}  
	\vspace{-0.03in}\\
	{\footnotesize (d) DCT-Gaussian-Sketch} & {\footnotesize (e) DCT-Gaussian-Sketch-PI} & {\footnotesize (f) rt-SVD}  \vspace{-0.03in}  \\
	{\footnotesize CPU: 0.08; PSNR: 13.03} & {\footnotesize CPU: 0.17; PSNR: 14.17} & {\footnotesize CPU: 0.21; PSNR: 13.68} \\ 
	
\end{tabular}
\caption{Qualitative results of the low-rank approximation performance of different methods on the HDU 2D code image with
	Gaussian white noise and the sketch size $k$ = 50 in terms of the CPU time and PSNR.}
\label{fig:LRA-WhiteNoise-2D-quali}	
\end{figure}

We now inspect the impact of noise on the performance of different methods. Speciﬁcally, the input tensor is randomly generated as
$$\mathcal{A} \leftarrow \mathcal{A}+\sigma \cdot \mathcal{E},$$
where $\mathcal{E}$ represents the Gaussian noise tensor following the standard Gaussian distribution, and $\sigma$ controls the noise level and is set to 5 here.

Figures \ref{fig:WhiteNoise-2D-decay-spec} and \ref{fig:LRA-WhiteNoise-2D-quali} respectively give the quantitative and qualitative results of the low-rank approximation performance of different methods on the HDU 2D code image (with size of $800\times800\times3$) with Gaussian white noise in terms of different metrics. It is again evident that as the sketch parameter size $k$ increases, the low-rank approximation performance of all the methods improves. Regarding the CPU time, our DCT-Gaussian-Sketch algorithm  outperforms all others. The DCT-Gaussian-Sketch algorithm also surpasses the T-Sketch algorithm in terms of accuracy, with reduced storage and manipulation requirements. The DCT-Gaussian-Sketch-PI algorithm achieves better accuracy by a large margin compared to the T-Sketch-PI algorithm with faster speed. Moreover, the DCT-Gaussian-Sketch-PI algorithm is only second to the truncated-t-SVD algorithm in accuracy, but it is remarkably faster than the truncated-t-SVD algorithm. In addition, from the right column of Figure \ref{fig:WhiteNoise-2D-decay-spec}, it can be seen that when $k\geq10$, the PSNR obtained by the DCT-Gaussian-Sketch and the DCT-Gaussian-Sketch-PI algorithms is higher than the original image with Gaussian white noise, indicating the effectiveness of our method both in tensor approximation and the denoising capacity.

\begin{figure}[H]
\centering
\includegraphics[trim={{.30\linewidth} {.08\linewidth} {.20\linewidth} {.08\linewidth}}, clip, width=0.87\linewidth, height = 0.47\linewidth]{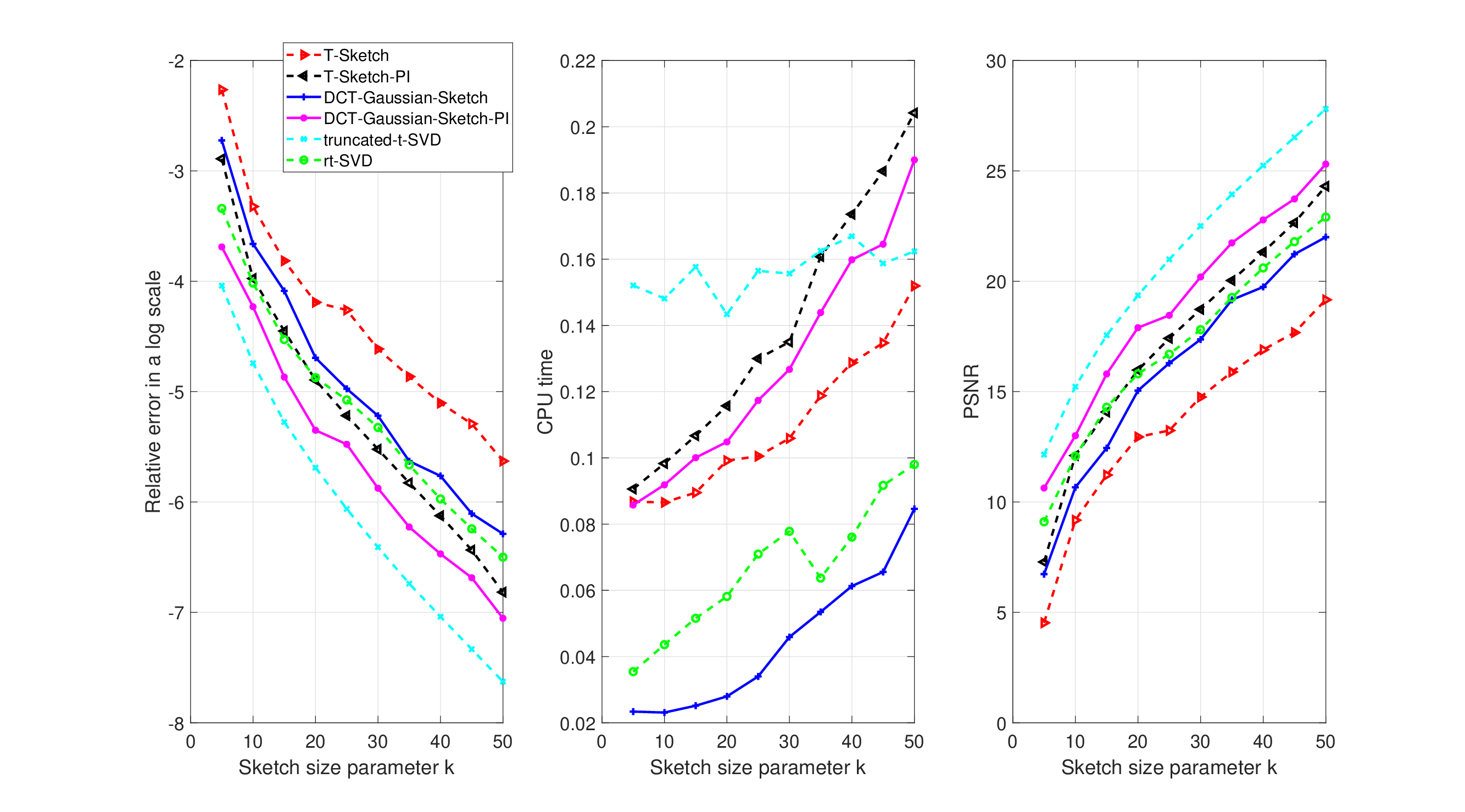}
\caption{Low-rank approximation performance of different methods on the `hall monitor' grayscale video clip (size: $144\times 176\times 30$) in terms of the relative error ({\it left} column), CPU time ({\it middle} column), and PSNR ({\it right} column).}
\label{fig:LRA-video-quan}
\end{figure}

\begin{figure}[htbp!]
\centering			       \begin{tabular}{ccc}
	\includegraphics[trim={{2.3in} {5.8in} {11.7in} {0.95in}}, clip, width=1.5in, height = 1.2in]{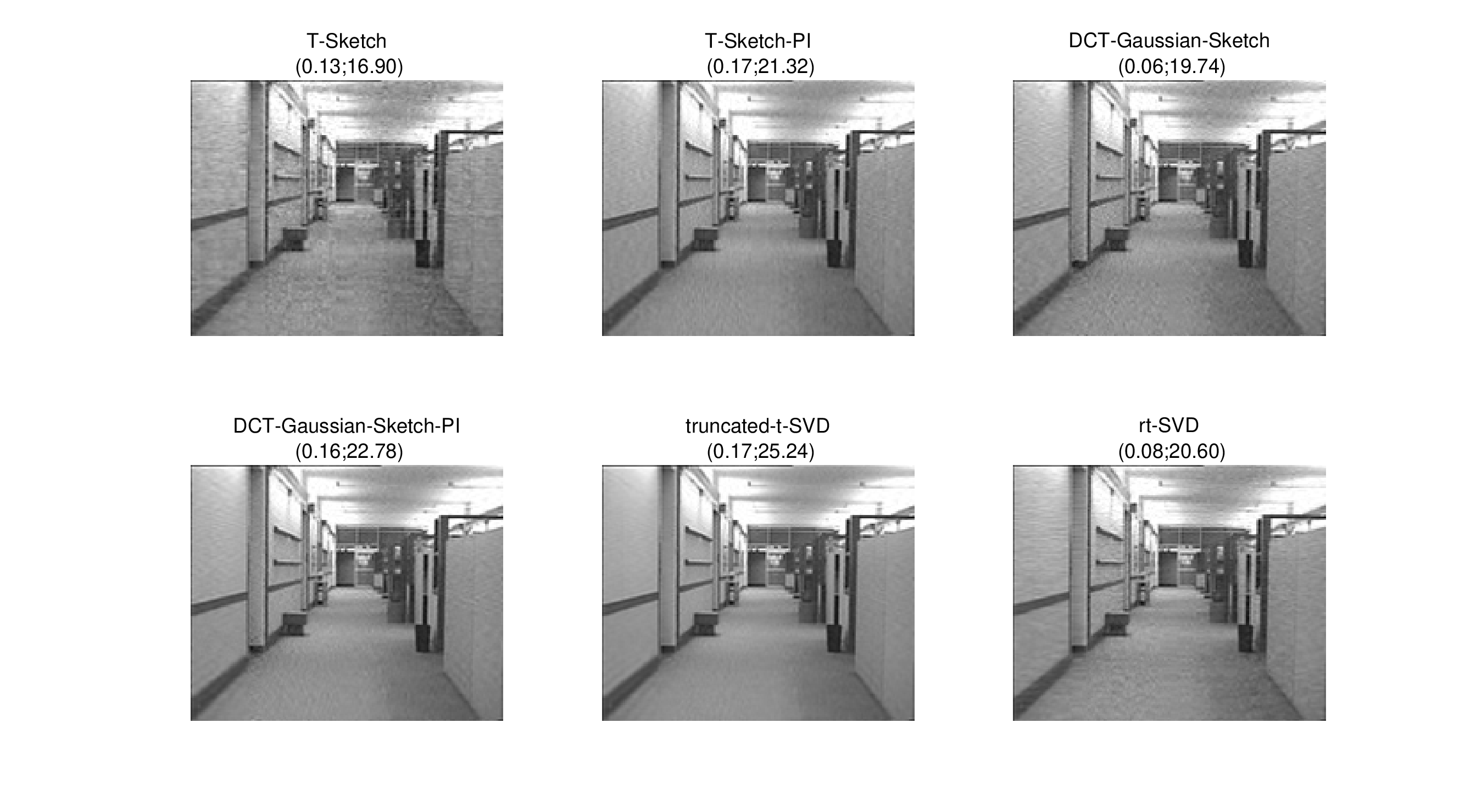} &
	\includegraphics[trim={{7.3in} {5.8in} {6.7in} {0.95in}}, clip, width=1.5in, height = 1.2in]{GaussianPIVideoK=40.eps} &
	\includegraphics[trim={{12.3in} {5.8in} {1.58in} {0.95in}}, clip, width=1.5in, height = 1.2in]{GaussianPIVideoK=40.eps} 
	\vspace{-0.03in}	\\
	{\footnotesize (a) T-Sketch} & {\footnotesize (b) T-Sketch-PI} & {\footnotesize (c) DCT-Gaussian-Sketch }  \vspace{-0.03in} \\
	{\footnotesize CPU:0.13; PSNR: 16.90}	& {\footnotesize CPU:0.17; PSNR: 21.32} & {\footnotesize CPU: 0.06; PSNR: 19.74} \\
	\includegraphics[trim={{2.3in} {1.1in} {11.7in} {5.62in}}, clip, width=1.5in, height = 1.2in]{GaussianPIVideoK=40.eps} &
	\includegraphics[trim={{7.3in} {1.1in} {6.7in} {5.62in}}, clip, width=1.5in, height = 1.2in]{GaussianPIVideoK=40.eps} &
	\includegraphics[trim={{12.3in} {1.1in} {1.58in} {5.62in}}, clip, width=1.5in, height = 1.2in]{GaussianPIVideoK=40.eps} 
	\vspace{-0.03in} \\
	{\footnotesize (d) DCT-Gaussian-Sketch-PI} & {\footnotesize (e) truncated-t-SVD} & {\footnotesize (f) rt-SVD}  \vspace{-0.03in}  \\
	{\footnotesize CPU: 0.16; PSNR: 22.78} & {\footnotesize CPU: 0.17; PSNR: 25.24} & {\footnotesize CPU: 0.08; PSNR: 20.60} 
\end{tabular}
\caption{Qualitative results of the low-rank approximation performance of different methods on the first frame of the `hall monitor' grayscale video clip with the sketch size $k$ = 40 in terms of the CPU time and PSNR.}
\label{fig:LRA-video-quali}	
\end{figure}

\subsubsection{Grayscale Video}	
We ultimately assess the proposed sketching algorithms and the peers using the extensively adopted YUV Video Sequences\footnote{http://trace.eas.asu.edu/yuv/index.html}. As an example, we utilize the `hall monitor' video clip and its first 30 frames to create a three-order tensor with dimensions $144\times176\times30$ for this experiment. 

Figures \ref{fig:LRA-video-quan} and \ref{fig:LRA-video-quali} respectively give the quantitative and qualitative results of the low-rank approximation performance of different methods on the `hall monitor' graysc-ale video clip in terms of different metrics. Consistent results are obtained as that obtained from the previous section for synthetic tensors and real-world images. For example, regarding the CPU time, our DCT-Gaussian-Sketch algorithm again outperforms all others.
Moreover, it also clearly surpasses the T-Sketch algorithm in terms of accuracy, with reduced storage and manipulation requirements. The DCT-Gaussian-Sketch-PI algorithm achieves better accuracy compared to the T-Sketch-PI algorithm with faster speed. 
The DCT-Gaussian-Sketch-PI algorithm is only second to the truncated-t-SVD algorithm in accuracy, but it is faster than the truncated-t-SVD algorithm. 

\section{Conclusion} \label{sec:con}
In this paper, we focused on large-scale tensor decomposition and proposed a novel two-sided sketching method based on the $\star{_{L}}$-product decomposition and transformed domains. Different transformed domains including the U, DCT, and DFT domains were investigated, together with different tensor sketching operators and the extension with the power iteration technique. A rigorous theoretical analysis was also conducted to assess the approximation error of the proposed method, particularly for the case of using the DCT transformation and the Gaussian tensor sketching operator with and without power iteration. Extensive numerical experiments and comparisons on low-rank approximation of synthetic large tensors and real-world data like color images and grayscale videos demonstrated the efficiency and effectiveness of the proposed approach in terms of both CPU time and low-rank approximation effects.

\section*{Appendix A}
This appendix is about random projection.
Three matrix sketching techniques, i.e., Gaussian projection, subsampled randomized Hadamard transform (SRHT), and count sketch \cite{{Wang15}}, are introduced below.

\subsection*{{A.1. Gaussian Projection}}
The Gaussian random projection matrix $S\in \mathbb{R}^{n \times s}$ is a matrix formed by $S = G/\sqrt{s}$, where
each entry of $G$ is sampled i.i.d. from $\mathcal{N} (0, 1)$. For matrix $A\in \mathbb{R}^{m \times n}$, the time complexity of Gaussian projection is $\mathcal{O}(mns)$. Algorithm \ref{alg:gp} presents the Gaussian projection process.

\begin{algorithm}{}
\caption{Gaussian Projection}
\begin{algorithmic}[1]
		\STATE   {\textbf{Input:}} $A \in \mathbb{R}^{m\times n}$, and parameter $s$. 
		\STATE {\textbf{function}} {\texttt{GaussianProjection}}{($A, s$)}
		\STATE Generate Gaussian random matrix $G \in \mathbb{R}^{n \times s}$;
		\STATE $S = \frac{1}{\sqrt{s}}G$;
		\STATE $C=AS$;
		\STATE {\textbf{return}} $C \in \mathbb{R}^{m\times s}$.
	\end{algorithmic}
	\label{alg:gp}
\end{algorithm}

\subsection*{A.2. Subsampled Randomized Hadamard Transform}
The SRHT matrix is defined by $S = DH_{n}P/\sqrt{sn} \in \mathbb{R}^{n \times s}$, where
\begin{itemize}
	\item $D \in \mathbb{R}^{n \times n}$
	is a diagonal matrix with diagonal entries sampled uniformly from $\{+1, 1\}$.
	\item $H_{n} \in \mathbb{R}^{n \times n}$
	is the Hadamard matrix defined recursively by $$H_{n}=\begin{pmatrix} H_{\frac{n}{2}} &H_{\frac{n}{2}} \\ H_{\frac{n}{2}} & -H_{\frac{n}{2}} \end{pmatrix}\quad {\rm and} \quad H_{2}=\begin{pmatrix} +1 &+1 \\ +1 & -1 \end{pmatrix} .$$ Note that the order of the Hadamard matrix is usually a power of 2.
	$\forall y \in \mathbb{R}^{n}$, the matrix vector product $y^{H}H_{n}$ can be performed in $\mathcal{O}({n \log n})$ by the fast Walsh-Hadamard transform algorithm in a divide-and-conquer fashion.
	\item $P \in \mathbb{R}^{n \times s}$ is a matrix that essentially samples $s$ columns from $DH_{n}$. 
\end{itemize}
Algorithm \ref{alg:srht} presents the SRHT process.


\begin{algorithm}{} 
	\caption{Subsampled Randomized Hadamard Transform (SRHT)}
	\begin{algorithmic}[1] 
		\STATE {\textbf{Input:}} $A \in \mathbb{R}^{m\times n}$, and parameter $s$.
		\STATE {\textbf{function}} \texttt{SRHT}{($A, s$)}
		\STATE Generate matrices $D \in \mathbb{R}^{n \times n}, H_n \in\mathbb{R}^{n \times n}$ and $P \in\mathbb{R}^{n \times s}$;
		\STATE $S = \frac{1}{\sqrt{sn}}DH_{n}P$;
		\STATE $C=AS$;
		\STATE {\textbf{return}} $C \in \mathbb{R}^{m\times s}$.
	\end{algorithmic}
	\label{alg:srht}
\end{algorithm}

\subsection*{A.3. Count Sketch}
There are two approaches, namely ``map-reduce fashion'' and ``streaming fashion'', to implement count sketch. Both approaches are equivalent, and we will focus on the streaming fashion here. The streaming fashion involves two steps: (i) the $m \times s$ matrix $C$ is initialized to zero; and (ii) for each column of $A$, its sign is randomly flipped with a probability of 0.5 and the flipped column is then added to a randomly chosen column of $C$. 

The count sketch process in the streaming fashion is detailed in Algorithm \ref{alg:cs-sf}. The streaming fashion keeps matrix $C$ in memory and processes the data $A$ in a single pass. When $A$ does not fit into memory, this approach is more efficient than the map-reduce fashion because it processes columns sequentially. Furthermore, if $A$ is a sparse matrix, sequentially accessing columns can be more efficient than random access.
It is noticed that the count sketch does not explicitly form the sketching matrix $S$ like the Gaussian projection and the SHRT processes. In fact, $S$ for count sketch is such a matrix that its each row has only one nonzero entry. The time complexity of count sketch is $\mathcal{O}(\text{nnz}(A))$, where $\text{nnz}(A)$ represents the number of nonzeros of matrix $A$.

\begin{algorithm}{}  \label{al3}
	\caption{Count Sketch in the Streaming Fashion}
	\begin{algorithmic}[1] 
		\STATE {\textbf{Input:}} $A \in \mathbb{R}^{m\times n}$, and parameter $s$.
		\STATE {\textbf{function}} \texttt{CountSketch}{($A,s$)}
		\STATE Initialize $C$ to be an $m \times s$ all-zero matrix;
		\STATE {\textbf{for}}{\ $i = 1 \;\textbf{to} \;n \ \textbf{do}$}
		\STATE \quad sample $l$ from the set $[s]$ uniformly at random;
		\STATE \quad sample $g$ from the set $\{+1, -1\}$ uniformly at random;
		\STATE \quad update the $l$-th column of $C$ by $C(:,l) \gets C(:,l) + g A(:,i)$; 
		\STATE {\textbf{end}}
		\STATE {\textbf{return}} $C \in \mathbb{R}^{m\times s}$. 
		
	\end{algorithmic}
	\label{alg:cs-sf}
\end{algorithm}

\section*{Appendix B}

In this appendix, 
we provide the proofs which are used to derive the error bound of the proposed two-sided sketching algorithms in Theorems \ref{thm4.7} and \ref{thmPI}.

\subsection*{B.1. Facts about Random Tensors} First, let us state a useful formula below that allows us to compute some expectations involving a Gaussian random tensor. 

\begin{proposition}
	\label{pro4.1}
	Assume $t > q+1$. Let $\mathcal{G}_{1}\in \mathbb{R}^{t\times q \times p}$ and $\mathcal{G}_{2}\in \mathbb{R}^{t\times l \times p}$ be Gaussian random tensors. For any tensor $\mathcal{B}$ with conforming dimensions,
	$$\mathbb{E}{\parallel \mathcal{G}_{1}^{\dagger} \star{_{L}} \mathcal{G}_{2} \star{_{L}} \mathcal{B} \parallel}_{F}^{2} = \frac{q}{t-q-1}{\parallel \mathcal{B} \parallel}_{F}^{2} \ .$$
\end{proposition}

\textbf{proof}\quad 
	By Lemma \ref{lem2.3} and the linearity of the expectation, we have
	\begin{eqnarray*}
		\mathbb{E}{\parallel \mathcal{G}_{1}^{\dagger} \star{_{L}} \mathcal{G}_{2} \star{_{L}} \mathcal{B} \parallel}_{F}^{2} = {1 \over p}\left(\sum_{i=1}^p \mathbb{E}\big\|\bar{\mathcal{G}_{1}^{\dagger}}^{(i)} \bar{\mathcal{G}_{2}}^{(i)} \bar{\mathcal{B}}^{(i)}\big\|_F^2\right).
	\end{eqnarray*}
	By using \cite[A.1]{Tropp19}, we have
	$$\mathbb{E}\big\|\bar{\mathcal{G}_{1}^{\dagger}}^{(i)} \bar{\mathcal{G}_{2}}^{(i)} \bar{\mathcal{B}}^{(i)}\big\|_F^2=\frac{q}{t-q-1}{\parallel \bar{\mathcal{B}}^{(i)} \parallel}_{F}^{2} \ ,$$
	which yields
	$$\mathbb{E}{\parallel \mathcal{G}_{1}^{\dagger} \star{_{L}} \mathcal{G}_{2} \star{_{L}} \mathcal{B} \parallel}_{F}^{2} = \frac{q}{t-q-1}{\parallel \mathcal{B} \parallel}_{F}^{2} \ .$$
	This completes the proof.

\subsection*{B.2. Results from Randomized Linear Algebra}
\begin{proposition}
	\label{pro4.2}
	Fix $\mathcal{A}\in \mathbb{R}^{m\times n\times p}$, let $\varrho < k$ be a natural number, and $f(\varrho,k) = {\varrho}/({k-\varrho-1})$.
	Then the tensor $\mathcal{Q}\in \mathbb{R}^{m\times k\times p}$ calculated by Eq. (\ref{eq3.12}) satisfies
	\begin{equation}
		\mathbb{E}_{\Omega}{\| \mathcal{A}-\mathcal{Q}\star{_{L}}\mathcal{Q}^{H}\star{_{L}}\mathcal{A} \|}_{F}^{2}\le (1+f(\varrho,k)){\tau^{2}_{\varrho+1}}(\mathcal{A}\star{_{L}}\mathcal{A}^{H}) \ .
	\end{equation}
	An analogous result holds for the tensor $\mathcal{P}\in \mathbb{R}^{n\times k\times p}$
	computed by Eq. (\ref{eq3.12}), i.e.,
	\begin{equation}
		\mathbb{E}_{\Upsilon}{\| \mathcal{A}-\mathcal{A}\star{_{L}} \mathcal{P}\star{_{L}} \mathcal{P}^{H} \|}_{F}^{2}\le (1+f(\varrho,k)){\tau^{2}_{\varrho+1}}(\mathcal{A}\star{_{L}}\mathcal{A}^{H}) \ .
	\end{equation}
\end{proposition}

\textbf{proof}\quad 
	By Lemma \ref{lem2.3} and the linearity of the expectation, we have
	\begin{eqnarray*}
		\mathbb{E}_{\Omega}\|\mathcal{A} - \mathcal{Q}_{L} \star{_{L}} \mathcal{Q}^{H}_{L} \star{_{L}} \mathcal{A} \|_F^2
		=  {1 \over p}\left(\sum_{i=1}^p \mathbb{E}\left\|\bar{\mathcal{A}}^{(i)}_{L} - \bar{\mathcal{Q}}^{(i)}_{L}\left(\bar{\mathcal{Q}}^{(i)}_{L}\right)^{H} \bar{\mathcal{A}}^{(i)}_{L}\right\|_F^2\right).
	\end{eqnarray*}
	By \cite[Theorem 10.5]{Halko11}, we have
	$$\mathbb{E}{\|{\bar{\mathcal{A}}}^{(i)}_{\text{$F_{n_{3}}$}}-\bar{\mathcal{Q}}^{(i)}_{\text{$F_{n_{3}}$}}(\bar{\mathcal{Q}}^{(i)}_{\text{$F_{n_{3}}$}})^{H}{\bar{\mathcal{A}}}^{(i)}_{\text{$F_{n_{3}}$}}\|}_{F}^{2}\leq (1+f(\varrho,k))\cdot \tau_{\varrho+1}^{2}(\bar{\mathcal{A}}^{(i)}_{\text{$F_{n_{3}}$}}(\bar{\mathcal{A}}^{(i)}_{\text{$F_{n_{3}}$}})^{H}) \ .$$
	Similarly, we know from \cite{Kernfeld15} that multiplying by a unitary matrix will not change the Frobenius norm, i.e., the Frobenius norm keeps the unitary matrix unchanged. Thus we have
	$$\mathbb{E}{\|{\bar{\mathcal{A}}}^{(i)}_{L}-\bar{\mathcal{Q}}^{(i)}_{L}(\bar{\mathcal{Q}}^{(i)}_{L})^{H}{\bar{\mathcal{A}}}^{(i)}_{L}\|}_{F}^{2}\leq (1+f(\varrho,k))\cdot \tau_{\varrho+1}^{2}(\bar{\mathcal{A}}^{(i)}_{L}(\bar{\mathcal{A}}^{(i)}_{L})^{H})$$
	and
	\begin{align*}
		\mathbb{E}_{\Omega}\|\mathcal{A} - \mathcal{Q}_{L} \star{_{L}} \mathcal{Q}^{H}_{L} \star{_{L}} \mathcal{A} \|_F^2
		& \leq \frac{1}{p}(1+f(\varrho,k))\cdot\big(\sum_{i=1}^{p}\tau_{\varrho+1}^{2}(\bar{\mathcal{A}}^{(i)}_{L}(\bar{\mathcal{A}}^{(i)}_{L})^{H})\big)\\
		& = (1+f(\varrho,k))\cdot \tau_{\varrho+1}^{2}(\mathcal{A}\star{_{L}} \mathcal{A}^{H}) \ .
	\end{align*}
	Similarly, there is 
	\begin{equation}
		\mathbb{E}_{\Upsilon}{\| \mathcal{A}-\mathcal{A}\star{_{L}} \mathcal{P}\star{_{L}} \mathcal{P}^{H} \|}_{F}^{2}
		\le 
		(1+f(\varrho,k)){\tau^{2}_{\varrho+1}}(\mathcal{A}\star{_{L}}\mathcal{A}^{H}) \ ,
	\end{equation}
	which completes the proof.

Let $\hat{\mathcal{A}}$ be the tubal rank $k$ approximation of $\mathcal{A}$ obtained by the L-Gaussian-Sketch algorithm. We now split the error ${\| \mathcal{A}-\hat{\mathcal{A}}\|} _{F}^{2}$ into two parts by Proposition \ref{pro4.3} below.
\begin{proposition} \label{pro4.3}
	Let $\hat{\mathcal{A}}$ be the tubal rank $k$ approximation of $\mathcal{A}\in \mathbb{R}^{m\times n\times p}$ obtained by the 
	L-Gaussian-Sketch algorithm, with $\mathcal{Q}$, $\mathcal{C}$ and $\mathcal{P}$ being the intermediate tensors obtained by the L-Gaussian-Sketch algorithm satisfying $\mathcal{Q}^{H}\star{_{L}} \mathcal{Q}=\mathcal{I}_{kkp}$ and $\mathcal{P}\star{_{L}} \mathcal{P}^{H}=\mathcal{I}_{nnp}$. Then
	\begin{equation}\label{eq4.21}
		{\| \mathcal{A}-\hat{\mathcal{A}}\|} _{F}^{2}={\| \mathcal{A}-\mathcal{Q}\star{_{L}}\mathcal{Q}^{H}\star{_{L}} \mathcal{A}\star{_{L}}\mathcal{P}\star{_{L}} \mathcal{P}^{H} \|}_{F}^{2}+{\| \mathcal{C}-\mathcal{Q}^{H}\star{_{L}} \mathcal{A}\star{_{L}} \mathcal{P} \|} _{F}^{2} \ .
	\end{equation}		
\end{proposition}

\textbf{proof}\quad
	Since $\hat{\mathcal{A}}=\mathcal{Q}\star{_{L}}\mathcal{C}\star{_{L}}\mathcal{P}^{H}$, we have
	\begin{eqnarray*}
		&& \|\mathcal{A} - \hat{\mathcal{A}} \|_F^2 \\
		& = &\frac{1}{p}\sum_{i=1}^{p}{\|\bar{\mathcal{A}}_{L}-\bar{\hat{\mathcal{A}}}_{L}\|}_{F}^{2} \\
		& = &\frac{1}{p}\sum_{i=1}^{p}({\|\bar{\mathcal{A}}_{L}-\bar{\mathcal{Q}}_{L}\bar{\mathcal{Q}}^{H}_{L}\bar{{\mathcal{A}}}_{L}\bar{\mathcal{P}}_{L}\bar{\mathcal{P}}^{H}\|}_{F}^{2}+
		{\|\bar{\mathcal{C}}_{L}-\bar{\mathcal{Q}}^{H}_{L}\bar{{\mathcal{A}}}_{L}\bar{\mathcal{P}}_{L}\|}_{F}^{2}) \\
		& = &{\| \mathcal{A}-\mathcal{Q}\star{_{L}}\mathcal{Q}^{H}\star{_{L}} \mathcal{A}\star{_{L}}\mathcal{P}\star{_{L}} \mathcal{P}^{H} \|}_{F}^{2}+{\| \mathcal{C}-\mathcal{Q}^{H}\star{_{L}} \mathcal{A}\star{_{L}} \mathcal{P} \|} _{F}^{2} \ .
	\end{eqnarray*}
	The first and last equations are known by Lemma \ref{lem2.3}, and the second equation is known by \cite[A.6]{Tropp19}. This completes the proof.

The error of the first part of Eq. \eqref{eq4.21}
can be given in the below Theorem \ref{thm4.4}.
\begin{theorem}\label{thm4.4}
	For any natural number $\varrho < k-1$, it holds that
	\begin{equation}
		\mathbb{E}_{\Upsilon}\mathbb{E}_{\Omega}{\| \mathcal{A}-\mathcal{Q}\star{_{L}} \mathcal{Q}^{H}\star{_{L}} \mathcal{A}\star{_{L}} \mathcal{P}\star{_{L}} \mathcal{P}^{H} \|}_{F}^{2}
		\le {(1+\frac{2\varrho}{k-\varrho-1})}\tau_{\varrho+1}^{2}(\mathcal{A}\star{_{L}}\mathcal{A}^{H}) \ .
	\end{equation}
\end{theorem}

\textbf{proof}\quad
	We have
	\begin{eqnarray*}
		&&\mathbb{E}_{\Upsilon}\mathbb{E}_{\Omega}{\| \mathcal{A}-\mathcal{Q}\star{_{L}} \mathcal{Q}^{H}\star{_{L}} \mathcal{A}\star{_{L}} \mathcal{P}\star{_{L}}\mathcal{P}^{H} \|}_{F}^{2}\\
		& = &\frac{1}{p}\sum_{i=1}^{p}\mathbb{E}_{\Upsilon}\mathbb{E}_{\Omega}{\|\bar{\mathcal{A}}_{L}-\bar{\mathcal{Q}}_{L}\bar{\mathcal{Q}}^{H}_{L}\bar{{\mathcal{A}}}_{L}\bar{\mathcal{P}}_{L}\bar{\mathcal{P}}^{H}\|}_{F}^{2}\\
		& \le &{(1+\frac{2\varrho}{k-\varrho-1})}\tau_{\varrho+1}^{2}(\mathcal{A}\star{_{L}}\mathcal{A}^{H}) \ .
	\end{eqnarray*}
	The first equation is known by Lemma \ref{lem2.3}, and the inequality is known by \cite[A.5]{Tropp19}. This completes the proof.

\subsection*{B.3. Decomposition of the Core Tensor Approximation Error}
We now obtain a formula for the error in the approximation $(\mathcal{C}-\mathcal{Q}^{H}\star{_{L}} \mathcal{A}\star{_{L}} \mathcal{P})$ (i.e., the second part of Eq. \eqref{eq4.21}).

Note that the core tensor $\mathcal{C}\in \mathbb{R}^{k\times k \times p}$
is defined in Eq. (\ref{eq3.15}), and the
orthonormal tensors $\mathcal{P}\in \mathbb{R}^{n\times k \times p}$ and $\mathcal{Q}\in \mathbb{R}^{m\times k \times p}$ are constructed in
Eq. (\ref{eq3.12}).
Let us introduce tensors $\mathcal{P}_{\perp}\in \mathbb{R}^{n\times (n-k)\times p}$ and $\mathcal{Q}_{\perp}\in \mathbb{R}^{m\times (m-k)\times p}$ whose ranges are complementary to those of $\mathcal{P}$ and $\mathcal{Q}$, respectively, i.e.,
\begin{align}	
	\mathcal{P}_{\perp}\star{_{L}} \mathcal{P}_{\perp}^{H} & = \mathcal{I}_{mmp}-\mathcal{P}\star{_{L}} \mathcal{P}^{H}, \\ \mathcal{Q}_{\perp}\star{_{L}} \mathcal{Q}_{\perp}^{H} & = \mathcal{I}_{mmp}-\mathcal{Q}\star{_{L}} \mathcal{Q}^{H},
\end{align}
where the columns of $\mathcal{P}_{\perp}$ and $\mathcal{Q}_{\perp}$ are orthonormal, separately. Next, we introduce the subtensors, i.e.,
\begin{equation}\label{eq4.22}
	\begin{split}
		\Phi_{1}:={\Phi}\star{_{L}} \mathcal{Q}\in \mathbb{R}^{s\times k\times p}, \ \ {\Phi}_{2}:=\Phi\star{_{L}} \mathcal{Q}_{\perp} \in \mathbb{R}^{s\times (m-k)\times p}, \\		{\Psi}_{1}^{H}:=\mathcal{P}^{H}\star{_{L}} {\Psi}^{H}\in \mathbb{R}^{k\times s\times p}, \ \ {\Psi}_{2}^{H}:={\mathcal{P}^{H}_{\perp}} \star{_{L}}{\displaystyle \Psi^{H}}\in \mathbb{R}^{(n-k)\times s\times p}.
	\end{split}
\end{equation}
With these notations at hand, we can state and prove Lemma \ref{lem4.5} below.
\begin{lemma}\label{lem4.5}
	(Decomposition of the Core Tensor Approximation) Assume the tub-
	al rank of $\Phi_{1}$ and ${\Psi}_{1}$ is $k$ and $s$, respectively. Then
	\begin{align*}
		& \ \mathcal{C}-\mathcal{Q}^{H}\star{_{L}} \mathcal{A}\star{_{L}} \mathcal{P}\\
		= & \ \Phi_{1}^{\dagger}\star{_{L}} \Phi_{2}\star{_{L}} (\mathcal{Q}^{H}_{\perp}\star{_{L}} \mathcal{A}\star{_{L}} \mathcal{P})+(\mathcal{Q}^{H}\star{_{L}} \mathcal{A}\star{_{L}} \mathcal{P}_{\perp})\star{_{L}} \Psi_{2}^{\dagger}\star{_{L}} ( \Psi^{\dagger}_{1})^{H}\\
		& +\Phi_{1}^{\dagger}\star{_{L}} \Phi_{2}\star{_{L}}(\mathcal{Q}^{H}_{\perp}\star{_{L}} \mathcal{A}\star{_{L}} \mathcal{P}_{\perp})\star{_{L}} \Psi_{2}^{\dagger}\star{_{L}} ( \Psi^{\dagger}_{1})^{H}.
	\end{align*}
\end{lemma}

\textbf{proof}\quad  
	Adding and subtracting terms, we write the core sketch $\mathcal{Z}$ as
	\begin{align*}
		\mathcal{Z}= & \ \Phi\star{_{L}} \mathcal{A}\star{_{L}} {\Psi}^{H}\\
		= & \ \Phi\star{_{L}} (\mathcal{A}-\mathcal{Q}\star{_{L}} \mathcal{Q}^{H}\star{_{L}} \mathcal{A}\star{_{L}} \mathcal{P}\star{_{L}} \mathcal{P}^{H})\star{_{L}} {\Psi}^{H}\\
		& + ( \Phi\star{_{L}} \mathcal{Q})\star{_{L}} (\mathcal{Q}^{H}\star{_{L}} \mathcal{A}\star{_{L}} \mathcal{P})\star{_{L}} (\mathcal{P}^{H}\star{_{L}} \Psi^{H}).
	\end{align*}
	Using Eq. (\ref{eq4.22}), we identify the tensors $\Phi_{1}$ and ${\Psi}_{1}$. For the above $\mathcal{Z}$, after left-multiplying it by $ \Phi^{\dagger}_{1}$ and right-multiplying it by $({ \Psi^{\dagger}_{1}})^{H}$, we have
	\begin{align*}
		\mathcal{C} & = \ {\Phi}^{\dagger}_{1}\star{_{L}} \mathcal{Z}\star{_{L}} ({ \Psi^{\dagger}_{1}})^{H}\\
		& = \ {\Phi}^{\dagger}_{1}\star{_{L}}  \Phi\star{_{L}} (\mathcal{A}-\mathcal{Q}\star{_{L}} \mathcal{Q}^{H}\star{_{L}} \mathcal{A}\star{_{L}} \mathcal{P}\star{_{L}} \mathcal{P}^{H})\star{_{L}}  \Psi^{H}\star{_{L}} ({ \Psi^{\dagger}_{1}})^{H} \\
		& \quad  \ \ +\mathcal{Q}^{H}\star{_{L}} \mathcal{A}\star{_{L}} \mathcal{P},
	\end{align*}  
	which identifies the core tensor $\mathcal{C}$ defined in (\ref{eq3.15}). For the above representation of $\mathcal{C}$, moving the term $\mathcal{Q}^{H}\star{_{L}}\mathcal{A}\star{_{L}}\mathcal{P}$ to the
	left-hand side will give the approximation error. Using  Eq. (\ref{eq4.22}) again, we have
	\begin{align*}	
		{\Phi}^{\dagger}_{1}\star{_{L}} \Phi & = {\Phi}^{\dagger}_{1}\star{_{L}}  \Phi\star{_{L}} \mathcal{Q}\star{_{L}} \mathcal{Q}^{H}+{\Phi}^{\dagger}_{1}\star{_{L}}  \Phi\star{_{L}} \mathcal{Q}_{\perp}\star{_{L}} \mathcal{Q}^{H}_{\perp} \\
		& = \mathcal{Q}^{H}+{ \Phi}^{\dagger}_{1}\star{_{L}} \Phi_{2}\star{_{L}}\mathcal{Q}^{H}_{\perp} \ ,
	\end{align*}
	and
	\begin{align*}
		{\Psi}^{H}\star{_{L}} ({\Psi}^{\dagger}_{1})^{H}
		= & \ \mathcal{P}\star{_{L}} \mathcal{P}^{H}\star{_{L}} {\Psi}^{H}\star{_{L}} ({ \Psi}^{\dagger}_{1})^{H}+\mathcal{P}_{\perp}\star{_{L}} \mathcal{P}_{\perp}^{H}\star{_{L}} {\Psi}^{H}\star{_{L}} ({\Psi}^{\dagger}_{1})^{H}\\
		= & \ \mathcal{P}+\mathcal{P}_{\perp}\star{_{L}}{{ \Psi}_{2}^{H}}\star{_{L}}({\Psi}^{\dagger}_{1})^{H}.
	\end{align*}
	Combining the last three displays, we have
	\begin{align*}
		& \ \mathcal{C}-\mathcal{Q}^{H}\star{_{L}} \mathcal{A}\star{_{L}} \mathcal{P}\\
		= & \ (\mathcal{Q}^{H}+{\Phi}^{\dagger}_{1}\star{_{L}} \Phi_{2}\star{_{L}} \mathcal{Q}_{\perp}^{H})\star{_{L}} (\mathcal{A}-\mathcal{Q}\star{_{L}} \mathcal{Q}^{H}\star{_{L}} \mathcal{A}\star{_{L}} \mathcal{P}\star{_{L}} \mathcal{P}^{H})\star{_{L}} (\mathcal{P}\\
		&+\mathcal{P}_{\perp}\star{_{L}} {{\Psi}_{2}^{H}}\star{_{L}}({ \Psi}^{\dagger}_{1})^{H}).
	\end{align*}
	Expanding the expression and using the orthogonality relations $\mathcal{Q}^{H}\star{_{L}} \mathcal{Q}=\mathcal{I}_{kkp}$, $\mathcal{Q}_{\perp}^{H}\star{_{L}} \mathcal{Q}=\mathcal{O}$, $\mathcal{P}^{H}\star{_{L}} \mathcal{P}=\mathcal{I}_{kkp}$, and $\mathcal{P}_{\perp}^{H}\star{_{L}} \mathcal{P}=\mathcal{O}$, we complete the proof.
	%

\subsection*{B.4. Probabilistic Analysis of the Core Tensor}	
We can then study
the probabilistic behavior of the error $(\mathcal{C}-\mathcal{Q}^{H}\star{_{L}} \mathcal{A}\star{_{L}} \mathcal{P})$,  conditional on $\mathcal{Q}$ and $\mathcal{P}$.

\begin{lemma}\label{lem4.6}
	(Probabilistic Analysis of the Core Tensor) Assume that the dimension reduction tensors $\Phi$ and $\Psi$ are Gaussian linear sketching operators. When $s \ge k$, it holds that
	\begin{equation}\label{eq4.23}
		\mathbb{E}_{\Phi,\Psi}[ \mathcal{C}-\mathcal{Q}^{H}\star{_{L}} \mathcal{A}\star{_{L}} \mathcal{P} ]=\mathcal{O}.
	\end{equation}
	When $s > k + 1$, the error can be expressed as
	\begin{align*}
		\mathbb{E}_{\Phi,\Psi}{\| \mathcal{C}-\mathcal{Q}^{T}\star{_{L}} \mathcal{A}\star{_{L}} \mathcal{P} \|} _{F}^{2}
		= & \ \frac{k}{s-k-1}{\| \mathcal{A}-\mathcal{Q}\star{_{L}} \mathcal{Q}^{H}\star{_{L}} \mathcal{A}\star{_{L}} \mathcal{P}\star{_{L}} \mathcal{P}^{H} \|} _{F}^{2} \\
		& \ +\frac{k(2k+1-s)}{(s-k-1)^{2}}{\| \mathcal{Q}_{\perp}^{H}\star{_{L}} \mathcal{A}\star{_{L}} \mathcal{P}_{\perp} \|} _{F}^{2} \ .
	\end{align*}
	In particular, when $s < 2k + 1$, the last term is nonnegative; and  when $s \ge 2k + 1$, the last term is nonpositive.
\end{lemma}	

\textbf{proof}\quad 
	Since $\Phi$ is a Gaussian linear sketching operator, the orthogonal subtensors $\Phi_{1}$ and $\Phi_{2}$ are also	Gaussian linear sketching operators because of the marginal property of the normal distribution. Likewise, $\Psi_{1}$ and $\Psi_{2}$ are Gaussian linear sketching operators. Provided that $s \ge k$, the tensor transformed tubal rank
	of $\Phi_{1}$ and $\Psi_{1}$ is $k$.
	Using the decomposition of the approximation error from Lemma \ref{lem4.5}, we have
	\begin{align*}
		\mathbb{E}_{\Phi,\Psi}[ \mathcal{C}-\mathcal{Q}^{H}\star{_{L}} \mathcal{A}\star{_{L}} \mathcal{P}]
		= & \ \mathbb{E}_{\Phi_{1}}\mathbb{E}_{\Phi_{2}}[ \Phi_{1}^{\dagger}\star{_{L}} \Phi_{2}\star{_{L}} (\mathcal{Q}_{\perp}^{H}\star{_{L}} \mathcal{A}\star{_{L}} \mathcal{P}) ]\\
		&\ +\mathbb{E}_{\Psi_{1}}\mathbb{E}_{\Psi_{2}}[ (\mathcal{Q}^{H}\star{_{L}} \mathcal{A}\star{_{L}} \mathcal{P}_{\perp})\star{_{L}} \Psi_{2}^{H}\star{_{L}}(\Psi_{1}^{\dagger})^{H} ]\\
		& \ + \mathbb{E}[ \Phi_{1}^{\dagger}\star{_{L}} \Phi_{2}\star{_{L}} (\mathcal{Q}_{\perp}^{H}\star{_{L}} \mathcal{A}\star{_{L}} \mathcal{P}_{\perp})\star{_{L}} \Psi_{2}^{H}\star{_{L}}( \Psi_{1}^{\dagger})^{H}] \ .
	\end{align*}
	Then we invoke independence to write the expectations as iterated expectations. Since $ \Phi_{2}$
	and $\Psi_{2}$ have mean zero. This formula makes it clear that the approximation error has 
	mean zero.
	
	To study the fluctuations, applying the independence and zero-mean property of $\Phi_{2}$
	and $\Psi_{2}$, we have
	\begin{align*}
		& \ \mathbb{E}_{\Phi,\Psi} {\| \mathcal{C}-\mathcal{Q}^{H}\star{_{L}} \mathcal{A}\star{_{L}} \mathcal{P} \|}_{F}^{2} \\
		= & \ \mathbb{E}_{\Phi} {\| \Phi_{1}^{\dagger}\star{_{L}} \Phi_{2}\star{_{L}} (\mathcal{Q}_{\perp}^{H}\star{_{L}} \mathcal{A}\star{_{L}} \mathcal{P}) \|}_{F}^{2} \\
		& \ +\mathbb{E}_{\Psi} {\| (\mathcal{Q}^{H}\star{_{L}} \mathcal{A}\star{_{L}} \mathcal{P}_{\perp})\star{_{L}} \Psi_{2}^{H}\star{_{L}} (\Psi_{1}^{\dagger})^{H} \|}_{F}^{2} \\
		& \ + \mathbb{E}_{\Phi}\mathbb{E}_{\Psi} {\| \Phi_{1}^{\dagger}\star{_{L}} \Phi_{2}\star{_{L}} (\mathcal{Q}_{\perp}^{H}\star{_{L}} \mathcal{A}\star{_{L}} \mathcal{P}_{\perp})\star{_{L}} \Psi_{2}^{H}\star{_{L}} (\Psi_{1}^{\dagger})^{H} \|}_{F}^{2} \ .
	\end{align*}
	Invoking Proposition \ref{pro4.1} yields
	\begin{align*}
		& \ \mathbb{E}_{\Phi,\Psi} {\| \mathcal{C}-\mathcal{Q}^{H}\star{_{L}} \mathcal{A}\star{_{L}} \mathcal{P} \|}_{F}^{2} \\
		= & \ \frac{k}{s-k-1}\Big[{\| \mathcal{Q}_{\perp}^{H}\star{_{L}} \mathcal{A}\star{_{L}} \mathcal{P} \|}_{F}^{2}+{\| \mathcal{Q}^{H}\star{_{L}} \mathcal{A}\star{_{L}} \mathcal{P}_{\perp} \|}_{F}^{2}\\
		&\ + \frac{k}{s-k-1}{\| \mathcal{Q}_{\perp}^{H}\star{_{L}} \mathcal{A}\star{_{L}} \mathcal{P}_{\perp} \|}_{F}^{2}\Big]\\
		= & \ \frac{k}{s-k-1}\Big[{\| \mathcal{Q}_{\perp}^{H}\star{_{L}} \mathcal{A}\star{_{L}} \mathcal{P} \|}_{F}^{2}+{\| \mathcal{Q}^{H}\star{_{L}} \mathcal{A}\star{_{L}} \mathcal{P}_{\perp} \|}_{F}^{2}+{\| \mathcal{Q}_{\perp}^{H}\star{_{L}} \mathcal{A}\star{_{L}} \mathcal{P}_{\perp} \|}_{F}^{2}\\
		& \ + \frac{2k+1-s}{s-k-1}{\| \mathcal{Q}_{\perp}^{H}\star{_{L}} \mathcal{A}\star{_{L}} \mathcal{P}_{\perp} \|}_{F}^{2}\Big] \ .
	\end{align*}
	Using the Pythagorean Theorem to combine the terms in the above equation, 
	we have
	\begin{align*}
		& \ \mathbb{E}_{\Phi,\Psi}{\| \mathcal{C}-\mathcal{Q}^{H}\star{_{L}} \mathcal{A}\star{_{L}} \mathcal{P} \|} _{F}^{2}\\
		= & \ \frac{k}{s-k-1}{\| \mathcal{A}-\mathcal{Q}\star{_{L}} \mathcal{Q}^{H}\star{_{L}} \mathcal{A}\star{_{L}} \mathcal{P}\star{_{L}} \mathcal{P}^{H} \|} _{F}^{2} \\
		& \ +\frac{k(2k+1-s)}{(s-k-1)^{2}}{\| \mathcal{Q}_{\perp}^{H}\star{_{L}} \mathcal{A}\star{_{L}} \mathcal{P}_{\perp} \|} _{F}^{2} \ .
	\end{align*}
	This completes the proof.

\end{document}